\documentclass[10pt]{article}
\usepackage{amsmath,amssymb,amsthm}
\usepackage{color}
\usepackage[nohead,margin=1.0in]{geometry}
\usepackage{booktabs}
\usepackage{threeparttable}
\usepackage{algorithm,algorithmic}
\usepackage{verbatim}
\usepackage{url}
\usepackage{graphicx}
\usepackage{epstopdf}
\usepackage{mathtools}
\usepackage{subfig}
\usepackage[colorlinks,linktocpage,linkcolor=blue]{hyperref}
\graphicspath{{./pics/}}
\allowdisplaybreaks

\usepackage{array}
\theoremstyle{plain}
\newtheorem{theorem}{Theorem}[section]
\newtheorem{remark}{Remark}[section]

\newtheorem{lemma}{Lemma}[section]

\newtheorem{corollary}{Corollary}[section]
\newtheorem{example}{Example}[section]

\numberwithin{equation}{section}

\def\al{\alpha}

\def\fy{\varphi}

\def\Om{\Omega}
\def\Dal{{\partial^\alpha_t}}

\def\S{\mathbb{S}}
\def\dH#1{\dot H^{#1}(\Omega)}

\def\A{\mathcal{A}}

\def\II{{(\Omega)}}

\renewcommand{\d}{\mathrm{d}}
\def\II{(\Omega)}

\title{Inverse Problems for Subdiffusion from Observation at an Unknown Terminal Time}
\author{Bangti Jin\thanks{Department of Mathematics, The Chinese University of Hong Kong, Shatin, New Territories, Hong Kong, P.R. China (\texttt{bangti.jin@gmail.com, btjin@math.cuhk.edu.hk}). The work of B. J. is supported by UK EPSRC grant EP/T000864/1 and EP/V026259/1, and a start-up fund from The Chinese University of Hong Kong.}
\and  Yavar Kian\thanks{Centre de Physique Th\'{e}orique  (CPT),  UMR-7332, Aix Marseille Universit\'{e}, Campus de Luminy, Case 907,
13288 Marseille cedex 9, France (yavar.kian@univ-amu.fr)}
\and Zhi Zhou\thanks{Department of Applied Mathematics, The Hong Kong Polytechnic University, Kowloon, Hong Kong, China. (\texttt{zhizhou@polyu.edu.hk}).
The work of Z. Zhou is partly supported by Hong Kong Research Grants Council (15303122) and an internal grant of Hong Kong Polytechnic University (Project ID: P0031041, Work Programme: ZZKS)}}

\begin{document}

\maketitle
\begin{abstract}
Inverse problems of recovering space-dependent parameters, e.g., initial condition, space-dependent source
or potential coefficient, in a subdiffusion model from the terminal observation have
been extensively studied in recent years. However, all existing studies have assumed that the terminal
time at which one takes the observation is exactly known. In this work, we present uniqueness and stability results
for three canonical inverse problems, e.g., backward problem, inverse  source and inverse potential
problems, from the terminal observation at an unknown time. The subdiffusive nature of the problem
indicates that one can simultaneously determine the terminal time and space-dependent parameter. The
analysis is based on explicit solution representations, asymptotic behavior of the
Mittag-Leffler function, and mild regularity conditions on the problem data. Further, we present several
one- and two-dimensional numerical experiments to illustrate the feasibility of the approach. 
\end{abstract}

\section{Introduction}

Let  $\Omega\subset\mathbb{R}^d $ ($d=1,2,3$) be an open bounded smooth domain with a boundary $\partial\Omega$.
Consider the following initial-boundary value problem with $\alpha\in(0,1)$ for the subdiffusion model:
 \begin{equation}\label{eqn:fde}
 \left\{\begin{aligned}
     \partial_t^\alpha u - \Delta u + q u&=f, &&\mbox{in } \Omega\times(0,\infty),\\
      u&=0,&&\mbox{on } \partial\Omega\times(0,\infty),\\
    u(0)&=u_0,&&\mbox{in }\Omega,
  \end{aligned}\right.
 \end{equation}
where $T>0$ is a fixed final time, $f\in L^\infty(0,T;L^2(\Omega))$ and $u_0\in L^2(\Omega)$ are given
source term and initial data, respectively, the nonnegative function $q\in L^\infty(\Omega)$ is a spatially
dependent potential, and $\Delta$ denotes the Laplace operator in space. The notation
$\Dal u(t)$ denotes the Djrbashian-Caputo fractional derivative in time $t$ of order
$\alpha\in(0,1)$ (\cite[p. 70]{KilbasSrivastavaTrujillo:2006} or \cite[Section 2.3.2]{Jin:2021book})
\begin{align}\label{eqn:RLderive}
   \Dal u(t)= \frac{1}{\Gamma(1-\alpha)}\int_0^t(t-s)^{-\alpha} u'(s)\d s,
\end{align}
where $\Gamma(z)$ is the Gamma function defined by $\Gamma(z) = \int_0^\infty s^{z-1}e^{-s}\d s$,
for $\Re (z)>0$. Note that the fractional derivative $\Dal u$ recovers the usual first-order derivative
$u'(t)$ as $\alpha\to1^-$ if $u$ is sufficiently smooth. Thus the model \eqref{eqn:fde} can be
viewed as a fractional analogue of the classical parabolic equation.

The model \eqref{eqn:fde} arises naturally in the study of anomalously slow diffusion processes,
which encompasses a broad range of important applications in engineering, physics and biology. The list of
successful applications include thermal diffusion in fractal media \cite{Nigmatulin:1986}, dispersion in
heterogeneous aquifer \cite{AdamsGelhar:1992}, ion dispersion in column experiments
\cite{HatanoHatano:1998}, and protein transport in membranes \cite{Kou:2008}, to name just a few. Thus its
mathematical theory has received immense attention in recent years; see the monographs
\cite{KubicaRyszewskaYamamoto:2020,Jin:2021book} for detailed discussions on the solution theory.
Related inverse problems have also been extensively studied \cite{JinRundell:2015,LiuLiYamamoto:2019review-source,
LiYamamoto:2019review}. The surveys \cite{LiuLiYamamoto:2019review-source} and \cite{LiYamamoto:2019review}
cover many inverse source problems and coefficient identification problems, respectively.

The observation $g(x)=u(x,T)$, $x\in \Omega$, at a terminal time $T$ is a popular choice
for the measurement data in practice. There is a very extensive literature on
inverse problems using terminal data, e.g., backward sudiffusion \cite{SakamotoYamamoto:2011,
ZhangZhou:2020}, inverse source problem \cite{JinRundell:2015,LiuLiYamamoto:2019review-source,JannoKinash:2018},
and inverse potential problem \cite{ZhangZhou:2017,KaltenbacherRundell:2019,JinZhou:2021ip,
ZhangZhangZhou:2022}, where the references are rather incomplete but we refer to the reviews
\cite{JinRundell:2015,LiuLiYamamoto:2019review-source,LiYamamoto:2019review} for further
references. Notably, several uniqueness and stability results have been proved. For example,
backward subdiffusion is only mildly ill-posed, and enjoys (conditional) Lipschitz stability \cite[Theorem
4.1]{SakamotoYamamoto:2011}, cf. \eqref{eqn:stab-u0} below. In all these existing studies, the
terminal time $T$ at which one collects the measurement has always been assumed to be fully known.
Nonetheless, in practice, the terminal time $T$ might be known only imprecisely. Therefore, it is
natural to ask whether one can still recover some information about the concerned parameter(s). The missing
knowledge of $T$ introduces additional challenges since the associated
forward map is not fully known then. In this work, we address this question in the affirmative both theoretically
and numerically, and study the inverse problem of identifying one of the following three parameters: (i)
initial condition $u_0$, (ii) space-dependent source component $\psi$, and (iii) space-dependent potential
$q$, from the observation $u(T)$ at an unknown terminal time $T$.

For each inverse problem, we shall establish the unique recovery of the space-dependent parameter
and the terminal time $T$ simultaneously from the terminal observation, as well as conditional stability estimates,
under suitable \textit{a priori} regularity assumptions on the initial data $u_0$
and the source $f$; see Theorems \ref{thm:stab-back}, \ref{thm:stab-source} and \ref{thm:stab-pot} for
the precise statements. The analysis relies heavily on explicit solution representations via Mittag-Leffler functions
(see, e.g., \cite{SakamotoYamamoto:2011}, \cite[Section 6.2]{Jin:2021book}). The essence of the argument is that the
regularity difference leads to distinct decay behavior of the Fourier coefficients of $u_0$ and $f$.
This combined with distinct polynomial decay
behaviour of Mittag--Leffler function $E_{\alpha,1}(z)$ (on the negative real axis) allows uniquely determining
the terminal time $T$. Note that the polynomial decay holds only for $E_{\alpha,1}(z)$ with a
order $\alpha\in (0,1)$, and does not hold in the integer case (i.e., $\alpha=1)$. Thus, the unique
determination of $T$ does not hold for normal diffusion. Once the terminal time
$T$ is determined, the unique determination of the space-dependent parameter follows. The
proof of the stability results relies on smoothing properties of the solution operators. In addition, we
present several numerical experiments to illustrate the feasibility of numerical recovery.
The numerical reconstructions are obtained using the Levenberg-Marquadt method
\cite{Levenberg:1944,Marquardt:1963}. Numerically, by choosing the hyperparameters in the method
properly, both space-dependent parameter and terminal time can be accurately recovered. To the  best of
our knowledge, this work presents the first uniqueness and stability results for inverse problems from
terminal data at an unknown time.

The rest of the paper is organized as follows. In Section \ref{sec:back}, we present uniqueness and stability
results for the backward problem, which are then extended to the inverse source problem in Section
\ref{sec:source}. In Section \ref{sec:q}, we discuss the inverse potential problem, which requires several
new technical estimates on the solution regularity and asymptotic decay. Finally,
some numerical results for one- and two-dimensional problems are given in Section
\ref{sec:numer}. Throughout, we denote by $u=u(v)$ and $\tilde u = u(\tilde v)$ the solutions to
problem \eqref{eqn:fde} with the space dependent parameter $v$ and $\tilde v$, respectively.
We often write a function $f(x,t):\Omega\times (0,T)\to \mathbb{R}$ as $f(t)$ a
vector-valued function on $(0,T)$. The notation $c$ denotes a generic constant which may differ at
each occurrence, but it is always independent of the concerned parameter and terminal time $T$.

\section{Backward problem}\label{sec:back}

In this section, we investigate the {backward problem} (\texttt{BP}): recover the initial data $u_0=u(0)$
from the solution profile $u(T)$ to problem \eqref{eqn:fde} at an unknown terminal time $T$.

\subsection{Solution representation}
First we recall the solution representation for problem \eqref{eqn:fde}, which plays a key role
in the analysis below. For any $s\ge0$, we denote by $\dH s\subset L^2(\Om)$ the Hilbert space
induced by the norm:
\begin{equation}\label{eqn:dH}
  \|v\|_{\dH s} =\Big(\sum_{j=1}^{\infty}\lambda_j^s ( v,\fy_j )^2\Big)^\frac12,
\end{equation}
with $\{\lambda_j\}_{j=1}^\infty$ and $\{\fy_j\}_{j=1}^\infty$ being respectively the eigenvalues (with multiplicity counted)
and eigenfunctions of the operator $A=-\Delta + qI$ on the domain $\Omega$
with a zero Dirichlet boundary condition. Then $\{\fy_j\}_{j=1}^\infty$ can be taken to form an orthonormal
basis in $L^2(\Omega)$. Further, $\|v\|_{\dH 0}$ is the norm in $L^2(\Omega)$,
$\|v\|_{\dH 1}$ is the norm in $H_0^1(\Om)$, and  $\|v\|_{\dH 2}=\|\Delta v\|_{L^2(\Omega)}$ is
equivalent to the norm in $H^2(\Om)\cap H^1_0(\Om)$ \cite[Section 3.1]{Thomee:2006}. For $s < 0$,
$\dot H^{s}\II$ denotes the dual space of $\dot H^{-s}\II$. Throughout, $(\cdot,\cdot)$ denotes both
duality pairing between $\dH {-s}$ and $\dH s$ and the $L^2(\Omega)$ inner product.

Now we represent the solution $u$ to problem \eqref{eqn:fde} using the eigenpairs $\{(\lambda_j,\fy_j)
\}_{j=1}^\infty$, following \cite{SakamotoYamamoto:2011} and \cite[Section 6.2]{Jin:2021book}.
Specifically, we define two solution operators $F(t)$ and $E(t)$ by
\begin{equation}\label{eqn:op}
F(t)v=\sum_{j=1}^\infty E_{\al,1}(-\lambda_jt^\al)(v,\fy_j)\fy_j\qquad \text{and}\qquad
E(t)v=\sum_{j=1}^\infty t^{\al-1}E_{\al,\al}(-\lambda_jt^\al)(v,\fy_j)\fy_j,
\end{equation}
where $E_{\alpha,\beta}(z)$ is the Mittag-Leffler function defined by
(\cite[Section 1.8, pp. 40-45]{KilbasSrivastavaTrujillo:2006} or \cite[Section 3.1]{Jin:2021book})
\begin{equation*}
  E_{\alpha,\beta}(z) = \sum_{k=0}^\infty \frac{z^{k}}{\Gamma(k\alpha+\beta)},\quad \forall z\in \mathbb{C}.
\end{equation*}
Then the solution $u$ of problem \eqref{eqn:fde} can be written as
\begin{equation}\label{eqn:sol-rep}
  u(t)=F(t)u_0+\int_0^t E(t-s)f(s)\d s.
\end{equation}

The function $E_{\alpha,\beta}(z)$ generalizes the familiar exponential function
$e^z$. The following decay estimates of $E_{\alpha,\beta}(z)$ are crucial in the analysis
below; See e.g., \cite[equation (1.8.28), p. 43]{KilbasSrivastavaTrujillo:2006} and
\cite[Theorem 3.2]{Jin:2021book} for the first estimate, and \cite[Theorem 4]{Simon:2014}
or \cite[Theorem 3.6]{Jin:2021book} for the second estimate.
\begin{lemma}\label{lem:ML-asymp}
Let $\alpha\in(0,2)$, $\beta\in \mathbb{R}$, and $\fy\in (\frac{\alpha\pi}{2},\min(\pi,\alpha\pi))$, and $N \in \mathbb{N}$.
Then for $\fy \le |\arg z| \le \pi$ with $|z| \rightarrow \infty$
$$ E_{\alpha,\beta}(z) = -\sum_{k=1}^N \frac{1}{\Gamma(\beta-\alpha k)} \frac1z + O\Big(\frac1{z^{N+1}}\Big). $$
For $0<\al_0<\alpha<\al_1<1,$ there exist constants $c_0$, $c_1>0$ depending only on $\al_0$ and $\al_1$ such that
\begin{align*}
 c_0({1-x})^{-1}\leq E_{\al,1}(x)\leq c_1({1-x})^{-1},\quad \forall x\leq 0.
\end{align*}
\end{lemma}

\subsection{Uniqueness and stability}

Now we study uniqueness and stability for \texttt{BP}, when the source
$f$ is time-independent, i.e., $f(x,t)\equiv f(x)$.
The key idea in proving uniqueness is to distinguish decay rates of Fourier coefficients of the initial
data $u_0$ (with respect to the eigenfunctions $\{\varphi_j\}_{j=1}^\infty$) and the
source $f$. We use the set $\S_\gamma$, $\gamma\in[-1,\infty)$, defined by
\begin{equation}\label{eqn:Sgamma}
   \S_\gamma =\Big \{ v\in \dH {-2}: ~\lim_{n\rightarrow \infty}\lambda_n^\gamma |(v,\fy_n)| = 0 \Big\}.
\end{equation}
Clearly, for any $\gamma \geq0$, $\dH \gamma \subset \S_{\frac{\gamma}{2}}$. If $v\in \dH {-2} \setminus \S_\gamma$, the
sequence $\{\lambda_n^\gamma |(v,\fy_n)|\}_{n=1}^\infty$ contains a subsequence that is bounded away from
zero, i.e., there exists $c_v>0$ and a sequence $\{n_\ell\}_{\ell=1}^\infty$ such that $\lim_{\ell\to\infty}n_\ell=\infty$ and $\lambda_{n_\ell}^\gamma|(v,\varphi_{n_\ell})|\geq c_v$, for all $\ell\in \mathbb{N}$.

\begin{theorem}\label{thm:back}
Let $f \in  \dH{-2}\, \backslash\, \S_\gamma$ for some $\gamma \ge 0$.  If \texttt{BP} has two solutions
$(T, u_{0})$ and $(\widetilde T,\widetilde u_{0})$ in the set  $\mathbb{R}_+\times \S_{\gamma+1}$ with the
observational data $u(T)$ and $\tilde u(\tilde T)$, respectively,  then $T= \widetilde T$ and $u_{0}=\widetilde u_{0}$.
\end{theorem}
\begin{proof}
Using the solution representation \eqref{eqn:sol-rep} and noting the identity
\begin{equation}\label{eqn:identity-1}
   \frac{\d}{\d t}E_{\alpha,1}(-\lambda_n t^{\alpha}) = -\lambda_nt^{\alpha-1}E_{\alpha,\alpha}(-\lambda_n t^{\alpha}),
\end{equation}
since $f$ is time-independent, the solution $u$ to problem \eqref{eqn:fde} can be written as
 \begin{equation*}
 \begin{aligned}
u(x,t) =   \sum_{n=1}^\infty \Big[E_{\alpha,1}(-\lambda_n  t ^{\alpha}) (u_0, \varphi_n )
+ \frac{1-E_{\alpha,1}(-\lambda_n  t ^\alpha)}{\lambda_n } (f,\varphi_n )\Big]  \varphi_n(x ).
  \end{aligned}
\end{equation*}
Let $\mathbb{K}=\{k\in\mathbb{N}:(f , \varphi_n)\neq0\}$, which under the condition
$f\in \dot H^{-2}(\Omega)\setminus\mathbb{S}_\gamma$ satisfies $|\mathbb{K}|=\infty$. For any $n \in \mathbb{K}$, taking
inner product (or duality pairing) with $\frac{\lambda_n\varphi_n}{(f,\varphi_n)}$ on both sides of
the identity gives
\begin{equation*}
 \begin{aligned}
\frac{\lambda_n (u(t)  , \varphi_n )}{(f, \varphi_n)}= \lambda_n E_{\alpha,1}(-\lambda_n  t ^{\alpha}) \frac{(u_0, \varphi_n )}{(f,\varphi_n )}
+ 1-E_{\alpha,1}(-\lambda_n  t ^\alpha) .
  \end{aligned}
\end{equation*}
Then setting $t=T$ and rearranging the terms lead to
\begin{equation}\label{eqb:relation-2}
 \begin{aligned}
\lambda_n \Big(1- \frac{\lambda_n (u(T)  , \varphi_n)}{(f,\varphi_n )}\Big) = - \lambda_n E_{\alpha,1}(-\lambda_n  T ^{\alpha}) \frac{\lambda_n(u_0, \varphi_n )}{(f,\varphi_n )}
+  \lambda_n E_{\alpha,1}(-\lambda_n  T ^\alpha).
  \end{aligned}
\end{equation}
By assumption, the source $f \in  \dH{-2}\, \backslash\, \S_\gamma$ and the initial condition $u_0 \in \S_{\gamma+1}$, and hence we deduce
\begin{equation}\label{eqn:limit-iv}
\lim_{n\in\mathbb{K},\rightarrow\infty} \frac{\lambda_n(u_0 , \varphi_n )}{(f,\varphi_n )}
= \lim_{n\in\mathbb{K},n\rightarrow\infty} \frac{\lambda_n^{\gamma+1}(u_0 , \varphi_n )}{\lambda_n^\gamma (f,\varphi_n )} =0.
\end{equation}
Then, by letting $n\rightarrow \infty$ and $t=T$, the relation \eqref{eqb:relation-2} implies
\begin{equation}\label{eqn:asymp-1-back}
 \begin{aligned}
\lim_{n\in \mathbb{K}, \ n\rightarrow \infty}\lambda_n \Big(1- \frac{\lambda_n (u(T)  , \varphi_n)}{(f,\varphi_n )}\Big)   =   \lim_{n\in \mathbb{K}, \ n\rightarrow \infty} \lambda_n E_{\alpha,1}(-\lambda_n  T^{\alpha})
= \frac{1}{\Gamma(1-\alpha)T^\alpha}.
  \end{aligned}
\end{equation}
The last identity follows from the fact that $\lambda_n(q)\rightarrow\infty$ and the asymptotic
behavior of $E_{\alpha,1}(z)$ in Lemma \ref{lem:ML-asymp}. Note that $\frac{1}{\Gamma(1-\alpha)
T^\alpha}$ is strictly decreasing in the terminal time $T$. Therefore, $T$ can be uniquely
determined from the observation $u(T)$. Finally, the unique determination of $u_0$ follows
from \cite[Theorem 4.1]{SakamotoYamamoto:2011}.
\end{proof}

\begin{remark}\label{rem:ass-1}
The validity of Theorem \ref{thm:back} relies crucially on the regularity difference
between the initial data $u_0$ and source $f$, so that the limit \eqref{eqn:limit-iv} holds.
\end{remark}

 \begin{remark}\label{rem:normal-1}
Theorem \ref{thm:back} shows the unique determination of the terminal time $T$ in problem
\eqref{eqn:fde} from the observation $u(T)$. This interesting phenomenon is due to
the distinct asymptotic behaviour of Mittag--Leffler functions and different smoothness of
the initial data $u_0$ and source  $f$. It sharply contrasts with the backward problem
of normal diffusion $(\alpha=1)$: analogous to \eqref{eqn:asymp-1-back}
\begin{equation*}
 \lim_{n\in \mathbb{K}, \ n\rightarrow \infty}\lambda_n \Big(1- \frac{\lambda_n (u(x,T)  , \varphi_n)}{(f,\varphi_n )}\Big)   =   \lim_{n\in \mathbb{K}, \ n\rightarrow \infty} \lambda_n e^{-\lambda_n  T}= 0 ,
\end{equation*}
and thus it is impossible to determine both the terminal time $T$ and initial
data $u_0$ from the data $u(T)$. This shows one distinct feature of anomalous slow
diffusion processes, when compared with standard diffusion.
\end{remark}

Next, we establish a stability estimate for \texttt{BP} with an approximately given $T$.
When the terminal time $T$ is exactly given, it recovers the following well-known estimate
\cite[Theorem 4.1]{SakamotoYamamoto:2011} (or \cite[Theorem 6.28]{Jin:2021book})
\begin{equation}\label{eqn:stab-u0}
  \|u_0-\tilde u_0\|_{L^2(\Omega)} \leq c\| u(T)- \tilde u(T)\|_{H^2(\Omega)}.
\end{equation}
\begin{theorem}\label{thm:stab-back}
Let $u_0$ and $\tilde u_0$ be the solutions of \texttt{BP} with observations $u(T)$ and $ \tilde u(\tilde T)$
with $T<\tilde T$, respectively. Then the following conditional stability estimate holds
\begin{equation*}
 \|u_0-  \tilde u_0 \|_{L^2\II} \le  c \big( 1+ c \tilde T^\alpha\big) \big(\| A(u(T)-u(\tilde T))  \|_{L^2\II}
+ c |T - \tilde T| T^{-1-\alpha}  \| u_0 - A^{-1} f \|_{L^2\II}\big).
\end{equation*}
\end{theorem}
\begin{proof}
By the solution representation \eqref{eqn:sol-rep} and the identity \eqref{eqn:identity-1}, we have
\begin{equation*}
 u(T) = F(T) u_0 + A^{-1}(I-F(T))f\quad
\text{and}\quad  \tilde u(\tilde T) = F(\tilde T) \tilde u_0 + A^{-1}(I-F(\tilde T))f.
\end{equation*}
Subtracting these two identities leads to
\begin{align*}
 u(T)  - \tilde u(\tilde T)
 &= \big(F(T) u_0 - F(\tilde T) \tilde u_0\big)+ A^{-1}\big(F(\tilde T)-F(T)\big)f \\
 &=F(\tilde T) (u_0 - \tilde u_0) + \big(F(T)- F(\tilde T)\big) (u_0-A^{-1}   f).
\end{align*}
Consequently, we arrive at
\begin{equation*}
   u_0 - \tilde u_0 = F(\tilde T)^{-1} \big(u(T)  - \tilde u(\tilde T) \big) - F(\tilde T)^{-1}  \big(F(T)- F(\tilde T)\big) (u_0-A^{-1}   f).
\end{equation*}
Next we bound the terms separately. For any $v\in L^2\II$, by Lemma \ref{lem:ML-asymp}, we derive
\begin{align*}
 \| A^{-1}F(\tilde T)^{-1} v \|_{L^2\II}^2 &= \sum_{n=1}^\infty \lambda_n^{-2}
 E_{\alpha,1}(-\lambda_n \tilde T^\alpha)^{-2}  (v, \fy_n)^2\nonumber \\
& \le  c \sum_{n=1}^\infty \Big(\frac{1+\lambda_n \tilde T^\alpha}{\lambda_n}\Big)^2  (v, \fy_n)^2
\le c\big( 1+ \tilde T^\alpha\big)^2 \| v \|_{L^2\II}^2 .
\end{align*}
Next we bound $F(\tilde T)^{-1}  (F(T)- F(\tilde T)) =   A^{-1}F(\tilde T)^{-1}  \int_{T}^{\tilde T} A F'(s) \,\d s$.
Note that for any $s\in (T,\tilde T)$, there holds \cite[Theorem 6.4]{Jin:2021book}
\begin{equation*}
   \|  A F'(s) v\|_{L^2\II} \le c s^{-1-\alpha} \|v \|_{L^2\II} \le c T^{-1-\alpha} \| v \|_{L^2\II}.
\end{equation*}
This estimate implies
\begin{equation}\label{eqn:AF-F'}
  \|A(F(T)- F(\tilde T))v\|_{L^2\II}  \leq c T^{-1-\alpha}|\tilde T-T| \| v \|_{L^2\II}.
\end{equation}
Consequently, we obtain
\begin{align*}
  \|F(\tilde T)^{-1}  (F(T)- F(\tilde T))v\|_{L^2\II}
   &\le   \|A^{-1}F(\tilde T)\| \|A(F(T)- F(\tilde T))v\|_{L^2\II} \\
  &\le  c \big( 1+ c\tilde T^\alpha\big) |T-\tilde T| T^{-1-\alpha} \| v \|_{L^2\II}.
\end{align*}
Then the desired result follows from immediately the preceding estimates.
\end{proof}

The next result bounds the terminal time $T$ in terms of data perturbation.
\begin{corollary}\label{cor:bwp}
Let $f \in  \dH{-2}\, \backslash\, \S_\gamma$ for some $\gamma \ge 0$. Let $(T, u_0), \, (\tilde T,\tilde u_0)
\in \mathbb{R}_+\times \S_{\gamma+1}$ be the solutions of \texttt{BP} with observations
$u(T)$ and $\tilde u(\tilde T)$, respectively. Then the following stability estimate holds
\begin{equation*}
 |T - \tilde T| \le c \min(\Lambda,\tilde\Lambda)^{-\frac1\alpha-1}|\Lambda -\tilde\Lambda|,
\end{equation*}
with the quantities $\Lambda$ and $\widetilde\Lambda$ respectively given by
\begin{equation*}
\Lambda = \lim_{n\in \mathbb{K}, \ n\rightarrow \infty}\lambda_n
\Big(1- \frac{\lambda_n (u(T)  , \varphi_n)}{(f,\varphi_n )}\Big)
\quad\mbox{and} \quad \tilde \Lambda = \lim_{n\in \mathbb{K}, \ n\rightarrow \infty}\lambda_n
\Big(1- \frac{\lambda_n (\tilde u(\tilde T)  , \varphi_n)}{(f,\varphi_n )}\Big).
\end{equation*}
In particular, for $\Lambda < \tilde \Lambda$, there holds
$\| u_0-  \tilde u_0 \|_{L^2\II} \le  c \big( 1+ c \tilde\Lambda^{-1} \big) \big(\| A(u(T)-\tilde u(\tilde T))  \|_{L^2\II}
  + c |\Lambda - \tilde \Lambda|  \|u_0 - A^{-1} f \|_{L^2\II}\big)$.
\end{corollary}
\begin{proof}
It follows from the relation \eqref{eqn:asymp-1-back} that
\begin{align*}
|T - \tilde T|  &= |\Lambda^{-\frac1\alpha} \Gamma(1-\alpha)^{-\frac1\alpha}
- \tilde \Lambda^{-\frac1\alpha} \Gamma(1-\alpha)^{-\frac1\alpha} |
\le \Gamma(1-\alpha)^{-\frac1\alpha}
\min(\Lambda,\tilde\Lambda)^{-\frac1\alpha-1}|\Lambda -\tilde\Lambda|.
\end{align*}
The assertion follows from Theorem \ref{thm:stab-back} and the identities
$T =  \Lambda^{-\frac1\alpha} \Gamma(1-\alpha)^{-\frac1\alpha}$
and  $\tilde T =  \tilde \Lambda^{-\frac1\alpha} \Gamma(1-\alpha)^{-\frac1\alpha}$.
\end{proof}

\section{Inverse source problem}\label{sec:source}

In this section, we extend the argument in Section \ref{sec:back} to an inverse source problem
of recovering the space dependent component from the observation $u(T)$. Following the standard
setup for inverse source problems \cite{LiuLiYamamoto:2019review-source}, we assume that the source
$f(x,t)$ is separable and satisfies
\begin{equation} \label{eqn:fg}
f(x,t) = g(t) \psi(x),\quad \text{with}~~ g\in L^\infty(0,T),~ g\ge c_g>0~~\text{and} ~~ \psi \in \dH{-1}.
\end{equation}
Then we consider the following inverse source problem (\texttt{ISP}): determine the spatially
dependent source component $\psi(x)$ from the solution profile $u(T)$ at a later but unknown time $T$.

First we give an intermediate result.
\begin{lemma}\label{lem:G-invertibility}
Let $G(T):=\int_0^{T} E(s) g(T-s)\,\d s$. Then under condition \eqref{eqn:fg}, $G$ is invertible and
\begin{equation*}
   \|  A^{-1}G(T)^{-1} v\|_{L^2\II} \leq c_g^{-1}{(1-E_{\alpha,1}(-\lambda_1T^\alpha))^{-1}}\|v\|_{L^2(\Omega)}.
\end{equation*}
\end{lemma}
\begin{proof}
Note that the function $E_{\al,\al}(-t) > 0$ for all $0\le t <\infty$ (since it is
completely monotone and analytic) (see, e.g., \cite{Schneider:1996}, \cite{MillerSamko:1998},
or \cite[Corollary 3.2 or Corollary 3.3]{Jin:2021book}).
This, the condition $ g\ge c_g >0$ and the identity \eqref{eqn:identity-1} imply
\begin{equation*}
  \int_0^{T}  s^{\alpha-1}E_{\alpha,\alpha}(-\lambda_ns^\alpha) g(T-s)\,\d s \geq c_g \int_0^T s^{\alpha-1}E_{\alpha,\alpha}(-\lambda_ns^\alpha)\d s
  = c_g\lambda_n^{-1}(1-E_{\alpha,1}(-\lambda_nT^\alpha)).
\end{equation*}
This implies the invertibility of the operator $G(T)$:
\begin{equation*}
  G(T)^{-1}v = \sum_{n=1}^\infty\frac{(v,\fy_n)}{
 \int_0^{T}  s^{\alpha-1}E_{\alpha,\alpha}(-\lambda_ns^\alpha) g(T-s)\,\d s}\varphi_n.
\end{equation*}
Consequently, for any $v\in L^2\II$, we have
\begin{align*}
 &\quad \|  A^{-1}G(T)^{-1} v \|_{L^2\II}^2 = \sum_{n=1}^\infty \Big[\frac{(v,\fy_n)}{\lambda_n
 \int_0^{T}  s^{\alpha-1}E_{\alpha,\alpha}(-\lambda_ns^\alpha) g(T-s)\,\d s}\Big]^2\\
 &\le  \sum_{n=1}^\infty \Big[\frac{(v,\fy_n)}{ c_g \lambda_n
 \int_0^{T}  s^{\alpha-1}E_{\alpha,\alpha}(-\lambda_ns^\alpha) \,\d s}\Big]^2
 = \sum_{n=1}^\infty \Big[\frac{(v,\fy_n)}{ c_g (1-E_{\alpha,1}(-\lambda_n T^\alpha))}\Big]^2\\
& \le \sum_{n=1}^\infty \Big[\frac{(v,\fy_n)}{(1-E_{\alpha,1}(-\lambda_1T^\alpha))c_g}\Big]^2 = c_g^{-2}(1-E_{\alpha,1}(-\lambda_1T^\alpha))^{-2}\|  v \|_{L^2\II}^2,
\end{align*}
where in the last inequality we have used the monotonicity of $E_{\alpha,1}(-t)\in (0,1]$
for $t\geq0$.
\end{proof}

The next result gives the unique determination of the source $\psi$ and the terminal time $T$.
\begin{theorem}\label{thm:source}
Let $u_0\in  L^2\II \, \backslash\, \S_\gamma$ for some $\gamma \ge 0$.  If \texttt{ISP} has two solutions
$(T, \psi)$ and $(\tilde T,\tilde \psi)$ in the set $\mathbb{R}_+\times \S_{\gamma}$ from the observations
$u(T)$ and $\tilde u(\tilde T)$, respectively, then $T= \tilde T$ and $\psi=\tilde \psi$.
\end{theorem}
\begin{proof}
Using the representation \eqref{eqn:sol-rep} and separability assumption \eqref{eqn:fg}, we have
 \begin{equation}\label{eqn:sol-repr-source}
u(t) =   \sum_{n=0}^\infty \Big[E_{\alpha,1}(-\lambda_n  t ^{\alpha}) (u_0, \varphi_n )
+ \int_0^t s^{\alpha-1} E_{\al,\al}(-\lambda_n  s^{\alpha})g(t-s)\,\d s\,(\psi,\varphi_n )\Big]  \varphi_n.
\end{equation}
Define $\mathbb{K}=\{k\in\mathbb{N}:\, (u_0, \varphi_n)\neq0\}$. For any $n \in \mathbb{K}$,
taking inner product (or duality pairing) with $\frac{\lambda_n \varphi_n}{(u_0,\varphi_n ) }$ on both
sides of the identity \eqref{eqn:sol-repr-source} and setting $t=T$, we have
\begin{align*}
\frac{  \lambda_n (u(T)  , \varphi_n )}{(u_0, \varphi_n)}=  \lambda_nE_{\alpha,1}(-\lambda_n  T^{\alpha})
+ \int_0^T s^{\alpha-1}E_{\al,\al}(-\lambda_n  s^{\alpha})g(T-s)\,\d s \frac{\lambda_n (\psi,\varphi_n)}{(u_0,\varphi_n)}.
\end{align*}
By assumption, $g\in L^\infty(0,T)$ and $E_{\al,\al}(-t) > 0$ for all $\infty>t\ge 0$ (since it is
completely monotone and analytic) \cite[Corollary 3.3]{Jin:2021book}, we have
\begin{equation*}
 \begin{aligned}
 \Big| \int_0^T s^{\alpha-1}E_{\al,\al}(-\lambda_n  s^{\alpha})g(T-s)\,\d s\Big|
  &\le   \int_0^T s^{\alpha-1}E_{\al,\al}(-\lambda_n  s^{\alpha}) \,\d s \|  g \|_{L^\infty(0,T)} \\
  & = \lambda_n^{-1}[1-E_{\alpha,1}(-\lambda_n T^{\alpha})] \|  g \|_{L^\infty(0,T)}
  \le \lambda_n^{-1}\|  g \|_{L^\infty(0,T)},
   \end{aligned}
\end{equation*}
where we have used the identity \eqref{eqn:identity-1} and the inequality
$E_{\alpha,1}(-t) \in [0,1]$ for all $t\ge0$. Since $\psi\in \S_\gamma$ and
$u_0\in L^2\II\backslash \S_\gamma$, we have
\begin{align*}
\lim_{n\in \mathbb{K}, \ n\rightarrow \infty} \Big|\int_0^T s^{\alpha-1}E_{\al,\al}(-\lambda_n  s^{\alpha})g(T-s)\,\d s
\frac{\lambda_n (\psi,\varphi_n)}{(u_0,\varphi_n)}\Big|
\le \lim_{n\in \mathbb{K}, \ n\rightarrow \infty}  \|  g \|_{L^\infty(0,T)} \Big|\frac{(\psi,\varphi_n)}{(u_0,\varphi_n)} \Big| = 0.
\end{align*}
Consequently, we have
\begin{equation}\label{eqn:asymp-1}
 \begin{aligned}
\lim_{n\in \mathbb{K}, \ n\rightarrow \infty}\frac{  \lambda_n (u(T)  , \varphi_n )}{(u_0, \varphi_n)}=
\lim_{n\in \mathbb{K}, \ n\rightarrow \infty} \lambda_nE_{\alpha,1}(-\lambda_n  T^{\alpha}) = \frac{1}{\Gamma(1-\alpha)T^\alpha}.
  \end{aligned}
\end{equation}
Note that the function $\frac{1}{\Gamma(1-\alpha)T^\alpha}$ is strictly decreasing in $T$. Hence, the
terminal time $T$ is uniquely determined by $u(T)$. Finally, the uniqueness of $\psi$
follows from the representation $\psi = G(T)^{-1}u(T)-G(T)^{-1}F(T)u_0$ (with $G(T)=\int_0^T E(s)
g(T-s)\,\d s$), and Lemma \ref{lem:G-invertibility}.
\end{proof}

\begin{remark}
Note that $\frac{1}{\Gamma(1-\alpha)T^\alpha}$ is also decreasing with respect to $\alpha$ if
$T$ is sufficiently large \cite[Lemma 4]{JannoKinash:2018}. So the terminal data $u(T)$ can uniquely
determine the fractional order $\alpha$ if $T$ is \textsl{a priori} known. See some related
argument in \cite{JannoKinash:2018, LiaoWei:2019} for the inverse source problem with an unknown
order $\alpha$.
\end{remark}

The next theorem gives a stability result for recovering the source $\psi$.
\begin{theorem}\label{thm:stab-source}
Fix $T_0>0$. Let $\psi$ and $\tilde \psi$ be the solutions of \texttt{ISP} with the data $u(T)$ and $\tilde u(\tilde T)$ with
$T_0\leq T<\tilde T$, respectively. Then for $g\in C^2[0,T]$, the following conditional stability estimate holds
\begin{equation*}
   \| \psi -  \tilde \psi \|_{L^2\II} \le  c  \big(\| A(u(T)-\tilde u(\tilde T))  \|_{L^2\II}
   +  |T -\tilde T| T^{-1-\alpha}  \| u_0\|_{L^2\II} +  |T - \tilde T| (1+ T^{-1} + T) \| \tilde \psi \|_{L^2\II} \big).
\end{equation*}
\end{theorem}
\begin{proof}
It follows from the solution representation \eqref{eqn:sol-rep} that
\begin{equation*}
  u(T) = F(T)u_0 + \int_0^{T} E(s) g(T-s)\,\d s \, \psi  \quad
 \text{and}\quad  \tilde  u(\tilde T) = F(\tilde T)u_0 + \int_0^{\tilde T} E(s) g(\tilde T-s)\,\d s \, \tilde  \psi.
\end{equation*}
Then subtracting these two identities leads to
\begin{align*}
 u(T)  - \tilde u(\tilde T) &= \big(F(T)  - F(\tilde T) \big) u_0 +  \int_0^{T} E(s) g(T-s)\,\d s \, (\psi  - \tilde \psi) \\
 &\quad + \int_0^{T} E(s) [g(T-s)- g(\tilde T-s)]\,\d s \tilde \psi  - \int_T^{\tilde T} E(s) g(\tilde T-s)\,\d s \tilde \psi.
\end{align*}
Therefore, with $G(T):=\int_0^{T} E(s) g(T-s)\,\d s$, we arrive at
\begin{align*}
\psi  - \tilde \psi & = G(T)^{-1} ( u(T)  - \tilde u(\tilde T) )
 +  G(T)^{-1}  \big(F(\tilde T) - F(T)  \big)u_0 \\
  &\quad + G(T)^{-1}\int_0^{T} E(s) [g(\tilde T-s)- g(T-s)]\,\d s \tilde \psi
+ G(T)^{-1} \int_T^{\tilde T} E(s) g(\tilde T-s)\,\d s \tilde \psi = \sum_{j=1}^4 {\rm I}_j.
\end{align*}
Next we bound the four terms separately. First, by Lemma \ref{lem:G-invertibility},
\begin{equation*}
 \|{\rm I}_1\|_{L^2\II} \le c_g^{-1}(1-E_{\alpha,1}(-\lambda_1T_0^\alpha))^{-1}\| A( u(T)  - \tilde u(\tilde T) )\|_{L^2\II}.
\end{equation*}
Second, by the estimate \eqref{eqn:AF-F'}, we can bound the term ${\rm I}_2$ by
\begin{equation*}
 \|{\rm I}_2  \|_{L^2\II} \le c  \| A  (F(\tilde T) - F(T)  ) u_0 \|_{L^2\II} \le c |\tilde T-T| T^{-1-\alpha} \| u_0\|_{L^2\II}.
\end{equation*}
Next we bound the term ${\rm I}_3$. It follows directly from integration by parts and the identities
$F'(t)=-AE(t)$ and $F(0)=I$ \cite[Lemmas 6.2 and 6.3]{Jin:2021book} that
\begin{align*}
 \int_0^{T} E(s) g(T-s)\,\d s &= A^{-1} \big(g(T)I - g(0) F(T)\big) -  A^{-1} \int_0^{T} F(s) g'(T-s)\,\d s,\\
\int_0^{T} E(s) g(\tilde T-s)\,\d s &= A^{-1} \big(g(\tilde T)I - g(\tilde T-T) F(T)\big)
-  A^{-1} \int_0^{T} F(s) g'(\tilde T-s)\,\d s.
\end{align*}
Since $g\in C^2[0,T]$, we have
\begin{align*}
  \|{\rm I}_3  \|_{L^2\II} &\le c \Big\|  A \int_0^{T} E(s) (g(T-s) - g(\tilde T-s))\,\d s \tilde \psi \Big\|_{L^2\II}\\
 & \le c \Big( |g(T) - g(\tilde T)| + |g(\tilde T - T) - g(0)| + \int_0^{T} |g'(\tilde T-s)  -  g'(T-s)|\,\d s\Big)
 \| \tilde \psi \|_{L^2\II}  \\
 &\le c |\tilde T-T| (1 + T) \|  g \|_{C^2[0,\tilde T]} \| \tilde \psi \|_{L^2\II}  .
\end{align*}
Finally, for the term ${\rm I}_4$, the estimate $\|AE(s)\|\leq cs^{-1}$ \cite[Theorem 6.4]{Jin:2021book} yields
\begin{align*}
 &\quad \|{\rm I}_4  \|_{L^2\II} \le c \Big\|   \int_{T}^{\tilde T} A E(s) g(\tilde T-s)\,\d s  \tilde \psi \Big\|_{L^2\II}\\
 &\le c   \|  g \|_{C[0,\tilde T]}  \int_{T}^{\tilde T} s^{-1} \,\d s \| \tilde \psi \|_{L^2\II}
 \le c   |\tilde T-T|  T^{-1} \|  g \|_{C[0,\tilde T]} \| \tilde \psi \|_{L^2\II}  .
\end{align*}
The preceding four estimates together complete the proof of the theorem.
\end{proof}

The next result bounds the terminal time $T$ for perturbed data. The proof is identical
with that for Corollary \ref{cor:bwp}, and hence it is omitted.
\begin{corollary}
Fix $T_0>0$. Let $u_0 \in L^2\II\, \backslash\, \S_\gamma$ for some $\gamma \ge 0$, and $(T, \psi), \,
(\tilde T,\tilde \psi) \in \mathbb{R}_+\times \S_{\gamma}$ with $T,\tilde T\geq T_0$ be the solutions of \texttt{ISP}
with observations $u(T)$ and $\tilde u(\tilde T)$, respectively. Then the following estimate holds
\begin{equation*}
 |T - \tilde T| \le  \Gamma(1-\alpha)^{-\frac1\alpha}
\min(\Lambda,\tilde\Lambda)^{-\frac1\alpha-1}|\Lambda -\tilde\Lambda|,
\end{equation*}
with the scalars $\Lambda$ and $\tilde \Lambda$ respectively given by
\begin{equation*}
\Lambda = \lim_{n\in \mathbb{K}, \ n\rightarrow \infty}
 \frac{\lambda_n (u(T)  , \varphi_n)}{(u_0,\varphi_n )}  \quad
\mbox{and}\quad \tilde \Lambda = \lim_{n\in \mathbb{K}, \ n\rightarrow \infty}
 \frac{\lambda_n (\tilde u(\tilde T)  , \varphi_n)}{(u_0,\varphi_n )} .
\end{equation*}
In particular, for $\Lambda < \tilde \Lambda$, there holds
$\| \psi -  \tilde \psi \|_{L^2\II} \le  c  \big(\| A(u(T)-\tilde u(\tilde T))  \|_{L^2\II}
  + |\Lambda -\tilde\Lambda| \Lambda^{-\frac1\alpha-1}(1+\Lambda^{\frac1\alpha}+\Lambda^{-\frac1\alpha}) \|  \tilde\psi \|_{L^2\II}\big)$.
\end{corollary}

\section{Inverse potential problem}\label{sec:q}
In this section, we discuss the identification of the potential $q$ in the
model \eqref{eqn:fde} from the observation $u(T)$, at an unknown terminal time $T$.
Specifically, we consider the domain $\Omega = (0,1)$ and a nonzero Dirichlet boundary condition:
\begin{equation}\label{eqn:fde-1}
  \left\{\begin{aligned}
     \partial_t^\alpha u - \partial_{xx} u +q u&=f, &&\mbox{in } \Omega\times(0,T],\\
     u(0,t)=a_0, u(1,t)&=a_1,&&\mbox{on } (0,T],\\
    u(0)&=u_0,&&\mbox{in }\Omega,
  \end{aligned}\right.
 \end{equation}
where the functions $f>0$ and $u_0>0$ are given spatially dependent source and
initial data, respectively, and $a_0$ and $a_1$ are positive constants. Throughout, the
potential $q$ belongs to the following admissible set
$$\A = \{ q\in L^\infty\II: \, 0 \le q\le c_0 \}.$$
The inverse potential problem (\texttt{IPP}) is to recover the
potential $q \in \mathcal{A}$ from the observation $u(T)$, for an unknown time $T$.

Similar to the discussions in Section \ref{sec:back}, let $A_q$ be the realization of the
elliptic operator $-\partial_{xx} + q(x) I$ in $L^2(\Omega)$, with its domain $\text{Dom}(A_q)$ given by
$\text{Dom}(A_q):=\{v\in L^2(\Omega): -\partial_{xx} v + qv\in L^2(\Omega)\ \text{and} \   v(0)= v(1)=0~~\text{in}~\partial\Omega\}.$
Let $\{(\lambda_{n}(q),\varphi_n(q))\}_{n=1}^\infty$ be the eigenpairs of the operator $A_q$, which is
not \textit{a priori} known for an unknown $q$. Note that for any $q\in \A$, the set $\{\varphi_n(q)\}_{n=1}^\infty$
can be chosen to form a complete (orthonormal) basis of the space $L^2\II$. It is well known that
for any $q\in\mathcal{A}$, the eigenvalues $\lambda_n(q)$ and eigenfunctions $\varphi_n(q)$
satisfy the following asymptotics \cite[Section 2, Chapter 1]{LevitanSargsjan:1975}:
\begin{equation}\label{eqn:asymp}
\lambda_n(q) = n^2\pi^2+O(1)\quad \text{and}\quad \fy_n(x;q) = \sqrt2\sin(n\pi x) + O(n^{-1}).
\end{equation}
Further, for any $v\in H_0^1\II\cap H^2\II$ and $q\in\A$,
the following two-sided inequality holds
\begin{equation}\label{eqn:equiv-n}
c_1\| v \|_{H^2\II} \le \| A_q v  \|_{L^2\II} + \| v  \|_{L^2\II} \le c_2\| v \|_{H^2\II},
\end{equation}
with constants $c_1$ and $c_2$ independent of $q$.
Next, we define an auxiliary function $\phi_q\in H^2(\Omega) $ satisfying
\begin{equation*}
\left\{\begin{aligned}
     &-\partial_{xx} \phi + q \phi = 0,\quad \text{in } \Omega=(0,1),\\
     & \phi_q(0) =a_0,\quad \phi_q(1) =a_1.
  \end{aligned}\right.
\end{equation*}
It is easy to see $\phi_0(x) = a_0(1-x)+a_1x$.
Then the solution $u$ of problem \eqref{eqn:fde-1} is given by
\begin{equation}\label{eqn:sol-rep-b}
   u(t) = F_q(t)u_0 + (I-F_q(t))\phi_q+(I-F_q(t)) A_q^{-1}f,
\end{equation}
where $E_q$ and $F_q$ denote the solution operators, cf. \eqref{eqn:op}, for
the elliptic operator $A_q$, and the subscript $q$ explicitly indicates the
dependence on the potential $q$.

Next, we show the unique recovery of the terminal time $T$. Like before, the key
is to distinguish the decay rates of Fourier coefficients of $u_0$ and $f$ with
respect to the (unknown) eigenfunctions $\{\varphi_j(q)\}_{j=1}^\infty$.
\begin{theorem}\label{thm:uniq-q}
Let $u_0\in L^2\II \backslash \dH {s'}$, with
$s'\in(0,\frac12)$, and $f\in \dH s$ and $q\in \A \cap \dH s$ with {$s \in( {s'} +\frac12,1)$}. Then in \texttt{IPP},
the terminal time $T$ is uniquely determined by the observation $u(T)$.
\end{theorem}
\begin{proof}
Let $\bar u_0=u_0-\phi_0$. 
Since $u_0 \in L^2 \II \backslash \dH {s'}$, by the asymptotics \eqref{eqn:asymp}, we obtain
\begin{equation}\label{eqn:Hs-bdd}
  \sum_{n=1}^\infty \lambda_n(q)^{s'}(\bar u_0, \sin(n\pi x))^2 = \infty.
\end{equation}
We claim that for {$\bar s = \frac{1}{4}+\frac{s'}{2}+\frac{\epsilon}{4}<\frac{s}{2}$}, with a small {$\epsilon\in(0,2s-2s'-1)$},
there holds $ \bar u_0 \notin \S_{\bar s}$. Indeed, assuming the contrary, i.e., $\bar u_0\in \S_{\bar s}$,
 the definition of $\S_{\bar s}$ in \eqref{eqn:Sgamma} implies 
$$\lim_{n\to \infty}\lambda_n(q)^{\bar s}|(\bar u_0,\sin(n\pi x))|=0,$$ 
and thus the sequence $\{\lambda_n(q)^{\bar s}|(\bar u_0,\sin(n\pi x))|\}_{n=1}^\infty$ 
is uniformly bounded. This and the asymptotics \eqref{eqn:asymp} lead to
\begin{equation*}
  \sum_{n=1}^\infty \lambda_n(q)^{s'}(\bar u_0, \sin(n\pi x))^2 = \sum_{n=1}^\infty \lambda_n(q)^{2\bar s} (\bar u_0, \sin(n\pi x))^2 \lambda_n(q)^{s'-2\bar s}
  \le c \sum_{n=1}^\infty  \lambda_n(q)^{s'-2\bar s}  \leq c\sum_{n=1}^\infty n^{-1-\epsilon}<\infty.
\end{equation*}
This contradicts the identity \eqref{eqn:Hs-bdd}, and hence the desired claim follows.
The claim $ \bar u_0 \notin \S_{\bar s}$ and the asymptotics \eqref{eqn:asymp} imply that there exists a constant
$c_*>0$, for any $N>0$, we can find $n>N$  such that
$n^{2\bar s}|(\bar u_0, \sin(n\pi x))| \ge c_*$. Let
$\mathbb{K}=\{n\in\mathbb{N}:~n^{2\bar s}|(\bar u_0, \sin(n\pi x))| \ge c_*\} $.
Then we have $|\mathbb{K}| = \infty$. By the asymptotics \eqref{eqn:asymp}, we may assume that
$ n^{2\bar s}|(\bar u_0, {\fy_n(q)})| \ge \sqrt{2}c_*/2$ for $n \in \mathbb{K}$, and
hence $(\bar u_0, {\fy_n(q)}) \neq 0$. 

Meanwhile, it follows directly from \eqref{eqn:sol-rep} that the solution $u(t)$ satisfies
\begin{equation}\label{eqn:exp-p}
u(t)-\phi_q =   \sum_{n=0}^\infty \Big[E_{\alpha,1}(-\lambda_n(q) t ^{\alpha}) (u_0-\phi_q , \varphi_n(q))
 + \big(1-E_{\alpha,1}(-\lambda_n(q) t ^\alpha)\big)
 \frac{(f,\varphi_n(q))}{\lambda_n(q)}\Big]   \varphi_n(q).
\end{equation} 
By taking inner product with $\varphi_n(q)$, $n\in\mathbb{K}$,
on both sides of the identity \eqref{eqn:exp-p}, we obtain
\begin{align}
   \frac{\lambda_n(q) (u(t)-\phi_0,\varphi_n(q))}{(\bar u_0, \varphi_n(q))}  &=  \lambda_n(q)E_{\alpha,1}(-\lambda_n(q) t ^{\alpha})  + \big(1-E_{\alpha,1}(-\lambda_n(q) t ^\alpha)\big)  \frac{(f,\varphi_n(q))}{(\bar u_0, \varphi_n(q))} \nonumber\\
    &\quad  +  \big(1-E_{\alpha,1}(-\lambda_n(q) t ^\alpha)\big)  \frac{\lambda_n(q) (\phi_q-\phi_0,\varphi_n(q))}{(\bar u_0, \varphi_n(q))}:=\sum_{i=1}^3{\rm I}_i.\label{eqn:relation-1}
\end{align}
We analyze the three terms ${\rm I}_i$, $i=1,2,3$, separately. Since $E_{\alpha,1}(-\lambda_nt^\alpha)\in(0,1]$
for $t\geq0$, the regularity condition $f\in \dH s$ implies $\lim_{n\to \infty}n^s |({f},  \varphi_n(q))|=0$. This and
the condition {$s>2\bar s$} yield
\begin{equation*}
{0\leq n^{2\bar s}|({f},  \varphi_n(q))| \leq n^{s}|({f},  \varphi_n(q))|}
\rightarrow0,\quad \mbox{as } n\to \infty.
\end{equation*}
Hence
\begin{align*}
\lim_{n\in\mathbb{K},\,n\to\infty}|{\rm I}_2| \le
  \lim_{n\in\mathbb{K},\,n\to\infty} \Big|\frac{ {n^{2\bar s}}(f,\varphi_n(q))}{ {n^{2\bar s}}(\bar u_0, \varphi_n(q))}  \Big|  = 0.
\end{align*}
By the definitions of the eigenpairs $(\lambda_n(q),\varphi_n(q))$ and
$\phi_q$ and $\phi_0$ and integration by parts, we have
\begin{equation*}
  \lambda_n(q) (\phi_q-\phi_0,\varphi_n(q)) = (\phi_q-\phi_0,-\partial_{xx}\varphi_n(q)+q\varphi_n(q))=-(q\phi_0,\varphi_n(q)).
\end{equation*}
Moreover, since $q \in \dH s$ with some {$s\in(\frac12,1)$} and
${\phi_0} \in C^\infty(\overline\Omega)$, we have $q\phi_0
\in \dH s$. Hence,
\begin{align}\label{eqn:bdd-I3}
\lim_{n\in\mathbb{K},\,n\to\infty}|{\rm I}_3| \le
 \lim_{n\in\mathbb{K},\,n\to\infty} \Big|  \frac{ ( q \phi_0,\varphi_n(q))}{(\bar u_0, \varphi_n(q))}   \Big|
 =\lim_{n\in\mathbb{K},\,n\to\infty} \Big|\frac{ {n^{2\bar s}}(q \phi_0,\varphi_n(q))}{ {n^{2\bar s}}(\bar u_0, \varphi_n(q))}  \Big|
  = 0 .
\end{align}
Now letting $n\rightarrow \infty$ and setting $t=T$, the relation \eqref{eqn:relation-1}
and the asymptotics of $E_{\alpha,1}(z)$ imply
\begin{equation}\label{eqn:relat-1}
 \begin{aligned}
\lim_{n\in \mathbb{K}, \ n\rightarrow \infty}\frac{\lambda_n(q) (u(T)-\phi_0,\varphi_n(q))}{(\bar u_0, \varphi_n(q))}  =   \lim_{n\in \mathbb{K}, \ n\rightarrow \infty} \lambda_n(q)E_{\alpha,1}(-\lambda_n(q) T^{\alpha}) = \frac{1}{\Gamma(1-\alpha)T^\alpha},
  \end{aligned}
\end{equation}
since $\lambda_n(q)\rightarrow\infty$, cf. \eqref{eqn:asymp}.
Note that the function $\frac{1}{\Gamma(1-\alpha)T^\alpha}$ is decreasing in
$T$ and independent of $q$. Last, we show that the left
hand side of \eqref{eqn:relat-1} can actually be computed independently of
$q \in \dH s \cap \A$, using the asymptotics \eqref{eqn:asymp}.
Indeed, since $u_0 \in L^2 \II \backslash \dH {s'}$,
$f\in \dH s$ and $q\in \dH s\cap \A$ with {$s \in ( {s'} +\frac12, 1)$ and  $s' \in (0, \frac12)$}, we have
$ A(q) u(T) \in L^2\II \backslash \dH {s'}$  \cite{SakamotoYamamoto:2011}.
This and the assumption $q\in \dH s$ imply $q u(T) \in \dH s$,
and hence $\partial_{xx} u(T) \in L^2\II \backslash \dH {s'}$.
Meanwhile, integration by parts twice yields
\begin{equation*}
  \lambda_n(q) (u(T)-\phi_0,\varphi_n(q)) = -(\partial_{xx}u(T),\varphi_n(q)) +(q(u(T) - \phi_0),\varphi_n(q)).
\end{equation*}
Noting the fact $q(u(T) - \phi_0)\in \dH s$ and then repeating the argument for \eqref{eqn:bdd-I3} yield
\begin{equation*}
\lim_{n\in \mathbb{K}, \ n\rightarrow \infty}\frac{(q(u(T) - \phi_0),\varphi_n(q))}{(\bar u_0 , \varphi_n(q))}=0.
\end{equation*}
Consequently, we derive
\begin{align}
\lim_{n\in \mathbb{K}, \ n\rightarrow \infty}\frac{\lambda_n(q) (u(T)-\phi_0,\varphi_n(q))}{(\bar u_0, \varphi_n(q))}
=\lim_{n\in \mathbb{K}, \ n\rightarrow \infty}\frac{  -(\partial_{xx}u(T),\varphi_n(q)) }{(\bar u_0, \varphi_n(q))}.\label{eqn:relat-2}
\end{align}
Now using the asymptotics \eqref{eqn:asymp} again, we obtain
\begin{align}
 \frac{  -(\partial_{xx}u(T),\varphi_n(q)) }{(\bar u_0, \varphi_n(q))}
 &=  \frac{  -(\partial_{xx}u(T), \sqrt{2}\sin(n\pi x)) + O(n^{-1}) \cdot \|  \partial_{xx}u(T) \|_{L^1\II} }{(\bar u_0, \sqrt{2}\sin(n\pi x))
 +  O(n^{-1}) \cdot \|\bar  u_0\|_{L^1\II}}\nonumber\\
 &= \frac{  - n^{{2\bar s} }(\partial_{xx}u(T), \sqrt{2}\sin(n\pi x)) + O(n^{{2\bar s-1}}) \cdot \|  \partial_{xx}u(T) \|_{L^1\II} }{n^{{2\bar s}}(\bar u_0, \sqrt{2}\sin(n\pi x))  + O(n^{{2\bar s}-1})\cdot \|\bar u_0\|_{L^1\II}}.\label{eqn:relat-3}
\end{align}
Since ${n^{2\bar s}}|(\bar u_0, \sin(n\pi x))| \ge c_*$ for all $n\in \mathbb{K}$, by the condition {$2\bar s-1<0$}, we obtain
\begin{equation*}
 \begin{aligned}
\lim_{n\in \mathbb{K}, \ n\rightarrow \infty}\frac{  -(\partial_{xx}u(T),\varphi_n(q)) }{(\bar u_0, \varphi_n(q))}
&=\lim_{n\in \mathbb{K} , \ n\rightarrow \infty}
\frac{- n^{{2\bar s}}(\partial_{xx}u(T), \sqrt{2}\sin(n\pi x)) + O(n^{{2\bar s-1}}) \cdot \|  \partial_{xx}u(T) \|_{L^1\II} }
 {n^{{2\bar s}}(\bar u_0, \sqrt{2}\sin(n\pi x))  + O(n^{{2\bar s-1}})\cdot \|\bar u_0\|_{L^1\II}}\\
 &=\lim_{n\in \mathbb{K} , \ n\rightarrow \infty}
\frac{- (\partial_{xx}u(T),  \sin(n\pi x))  }{ (\bar u_0,  \sin(n\pi x))}.
  \end{aligned}
\end{equation*}
Therefore, the terminal time $T$
is uniquely determined by the observation $u(T)$.
\end{proof}
\begin{remark}
The independence of the limit in \eqref{eqn:relat-1} on the potential
$q$ relies on the asymptotics \eqref{eqn:asymp}. This seems valid only in the one-dimensional
case, and it represents the main obstacle for the extension to the multi-dimensional case.
Nonetheless, the rest of the analysis does not use the estimate \eqref{eqn:asymp}, and all
the remaining results hold also for the multi-dimensional case.
\end{remark}

Next we determine the potential $q\in \A$ from the observation $u(T)$. First, we give useful
smoothing properties of the operators $F_q$ and $E_q$. The notation $\|\cdot\|$ denotes the
operator norm on $L^2(\Omega)$.
\begin{lemma}\label{lem:op}
For $q\in\A$, there exists a $c>0$ independent of $q$ and $t$ such that for any $s=0, 1$ and $\ell=0, 1$,
\begin{align*}
 \|F_q(t)\|+ t^{1-\alpha}\|E_q(t)\| \le c \min(1,   t^{-\alpha})\quad \mbox{and}\quad
   \|A_q^s   F_q^{(\ell)}(t)\| \le c t^{-\alpha-\ell}.
\end{align*}
\end{lemma}
\begin{proof}
The estimates follow from \cite[Theorems 6.4]{Jin:2021book}. Indeed, by \cite[Theorem 6.4(iv)]{Jin:2021book}
and Lemma \ref{lem:ML-asymp},
\begin{equation*}
  \|E_q(t)\|\le t^{\alpha-1}E_{\alpha,\alpha}(-\lambda_1(q)t^\alpha) \leq ct^{\alpha-1}\min(1, t^{-\alpha}).
\end{equation*}
The bound on $F_q(t)$ follows similarly. For the second estimate, the case $s=1$, $\ell=0,1$, the assertion is
contained in \cite[Theorem 6.4(iii)]{Jin:2021book}, and the case $s=0, \ell=0$ is direct from the first estimate.
The remaining case $s=0,\ell=1$ follows from Lemma \ref{lem:ML-asymp} (noting $1/\Gamma(0)=0)$:
\begin{align*}
  \|F_q'(t)v\|_{L^2\II}^2 & = \|A_qE_q(t)v\|_{L^2(\Omega)}^2 = \sum_{n=1}^\infty \lambda_n(q)^2t^{2\alpha-2} E_{\alpha,\alpha}(-\lambda_n(q)t^\alpha)^2(v,\varphi_n(q))^2 \\
  &\leq ct^{2\alpha-2}\sum_{n=1}^\infty \frac{\lambda_n^2}{(1+(\lambda_nt^\alpha)^2)^2}(v,\varphi_n(q))^2 \leq ct^{-2-2\alpha}\|v\|_{L^2(\Omega)}^2.
\end{align*}
Combining these assertions completes the proof of the lemma.
\end{proof}

The next lemma gives \textit{a priori} estimate on the solution $u$ to problem \eqref{eqn:fde-1}.
\begin{lemma}\label{lem:Dalu}
Let $u_0, f\in L^2\II$ and $q\in \A$, and $u$ be the solution to problem \eqref{eqn:fde-1}. Then
there exists $c>0$ independent of $q$ and $t$ such that
\begin{align*}
  \|  \partial_t u(t) \|_{H^2\II} &\le c t^{-\alpha-1},\\
 \| \Dal u  (t)\|_{L^2\II} + \|  u(t) \|_{H^{2}\II} &\le c (1 + t^{-\alpha}).
\end{align*}
\end{lemma}
\begin{proof}
The proof employs \eqref{eqn:equiv-n} and Lemma \ref{lem:op}. First, by
\eqref{eqn:sol-rep-b}, we have $\partial_t u(t)   = F_q'(t)(u_0- \phi_q- A_q^{-1}f).$ Then from
Lemma \ref{lem:op} and the norm equivalence \eqref{eqn:equiv-n}, the first estimate follows
\begin{align*}
\| \partial_t u(t) \|_{H^2 \II} &\le c(\|  A_q F_q'(t)\|+\|F'_q(t)\|)\|u_0 -\phi_q- A_q^{-1} f\|_{L^2\II} \le c t^{-\alpha-1}.
\end{align*}
Similarly, by \eqref{eqn:sol-rep-b} and Lemma \ref{lem:op}, there holds
\begin{align*}
\| u(t) \|_{L^2\II}  &\le \| F_q(t)u_0 + (I-F_q(t))\phi_q+(I-F_q(t)) A_q^{-1}f \|_{L^2\II}\leq c,\\
\| u(t) \|_{H^2\II}  &\le \| F_q(t) (u_0-\phi_q)+(I-F_q(t)) A_q^{-1}f \|_{\dot H^2\II} + \|\phi_q  \|_{H^2\II}\\
&\le c \min(1,t^{-\alpha}) \|  u_0-\phi_q \|_{L^2\II} + \| f \|_{L^2\II} + \|\phi_q \|_{H^2\II} \le c (1 + t^{-\alpha}).
\end{align*}
The bound on $\partial_t^\alpha u$ is direct from the identity
$\Dal  u(t) =  -A_q F_q(t) (u_0 +  \phi_q  -  A_q^{-1}f)$
and Lemma \ref{lem:op}.
\end{proof}

For any $q\in \A$, we denote the solution $u$ to problem \eqref{eqn:fde-1} by $u(q)$.
The next lemma provides a crucial \textit{a priori} estimate.
Like before, we denote by $u$ and $\tilde u$ to be $u(q)$ and $u(\tilde q)$ below.
\begin{lemma}\label{lem:stab-0}
Let $T_*>0$ be fixed, and $q, \tilde q\in \A$. Then there exists $c>0$ independent of $q$, $\tilde q$ and $t$ such that for any $t\ge T_* > 0$,
\begin{equation*}
\|\Dal(u-\tilde u)(t)\|_{H^2\II} \le ct^{-\alpha} \|q-\tilde q\|_{L^2\II}.
\end{equation*}
\end{lemma}
\begin{proof}
Let $w=u - \tilde u$. Then $w$ solves
 \begin{equation}\label{PDE-w}
 \left\{ \begin{aligned}
    \Dal w-\Delta w+q w&=(\tilde q-q)\tilde u, &&\mbox{in }\Omega\times(0,T],\\
    w&=0,&&\mbox{on }\partial\Omega\times(0,T],\\
    w(0)&=0,&&\mbox{in }\Omega.
  \end{aligned}\right.
 \end{equation}
The representation \eqref{eqn:sol-rep} implies $w(t) = \int_0^t E_q(t-s) (\tilde q-q)\tilde u(s) \,\d s$.
The governing equation for $w$ and the identities $\Dal F_q(t) = -A_q F_q(t)$ and $A_qE_q(t) = - F_q'(t)$
\cite[Lemma 6.3]{Jin:2021book} lead to
\begin{equation*}
\begin{split}
\partial_t^\alpha w(t) &=  -  A_q \int_0^t E_q(t-s) (\tilde q-q)\tilde u(s)  \,\d s +   (\tilde q-q)\tilde u(t) \\
&= \int_0^t F_q'(t-s) (\tilde q-q)\tilde u(s)  \,\d s +   (\tilde q-q)\tilde u(t)  = \partial_t \int_0^t F_q(t-s) (\tilde q-q)\tilde u(s) \,\d s.
\end{split}
\end{equation*}
Let $\phi(t) =  \int_0^t F_q(t-s) (\tilde q-q)\tilde u(s) \,\d s$. Sobolev embedding theorem
and Lemmas \ref{lem:op} and \ref{lem:Dalu} imply
\begin{align*}
\|A_q \phi(t) \|_{L^2\II}
&\le  \| \tilde q-q \|_{L^2\II}\int_0^t \| A_q F_q(t-s)  \| \ \| \tilde u(s) \|_{L^\infty\II}\, \d s\\
&\le  c \| \tilde q-q \|_{L^2\II}\int_0^t (t-s)^{-\alpha}  \|\tilde u(s) \|_{H^{2 }\II}\, \d s \\
&\le c (t^{1-\alpha} + t^{1-2\alpha}) \| \tilde q-q \|_{L^2\II}
\le c_{T_*}  t^{1-\alpha}  \| \tilde q-q\|_{L^2\II} .
\end{align*}
Next, the identity
$t A_q \phi(t) = \int_0^t (t-s)A_q  F_q(t-s) (\tilde q-q)\tilde u(s) \,\d s  + \int_0^t A_qF_q(s) (\tilde q-q) (t-s)\tilde u(t-s) \,\d s$ yields
\begin{align*}
 \partial_t(t  A(q_1)  \phi(t))  &= \int_0^t  \big[(t-s) A_q F_q'(t-s) +  A_q F_q(t-s)\big] (\tilde q-q)\tilde u(s) \,\d s \\
 &\quad + \int_0^t  A_q  F_q(s) (\tilde q-q)  \big[\tilde u(t-s)+ (t-s)\tilde u'(t-s)\big] \,\d s =: {\rm I}_1 + {\rm I}_2.
\end{align*}
Next we derive bounds for ${\rm I}_1$ and ${\rm I}_2$. First, we bound the term ${\rm I}_1$ by
Lemmas \ref{lem:op}  and \ref{lem:Dalu}:
\begin{align*}
  \| {\rm I}_1 \|_{L^2\II} &\le \int_0^t \big[(t-s)\|A_q F_q'(t-s)\| + \|A_q F_q(t-s)\| \big]
  \|(\tilde q-q)\tilde u(s)\|_{L^2\II} \,\d s\\
  &\le c\| \tilde q-q \|_{L^2\II}  \int_0^t  (t-s)^{-\alpha}  \|\tilde u(s)\|_{L^\infty\II} \,\d s
  \le c_{T_*}  t^{1-\alpha}  \| \tilde q-q \|_{L^2\II} .
\end{align*}
Similarly, by Lemma \ref{lem:op}, Sobolev embedding theorem and Lemma \ref{lem:Dalu}, the term ${\rm I}_2$ can be bounded as
\begin{align*}
 \|{\rm I}_2 \|_{L^2\II} &\le \int_0^t  \|A_qF_q(s)\| \, \|\tilde q-q\|_{L^2\II}  \big[(t-s)\|\tilde u'(t-s)\|_{L^\infty\II} +\|\tilde u(t-s)\|_{L^\infty\II}\big]  \,\d s\\
  &\le c \| \tilde q-q \|_{L^2\II} \int_0^t  s^{-\alpha} \big[(t-s)\|\tilde u'(t-s)\|_{H^{2}\II} +\|\tilde u(t-s)\|_{H^{2}\II}\big]    \,\d s\\
  &\le c   \| \tilde q-q\|_{L^2\II}  \int_0^t s^{-\alpha} ((t-s)^{-\alpha} + 1)  \,\d s\le c_{T_*}  t^{1-\alpha}  \|\tilde q-q\|_{L^2\II} .
\end{align*}
Then the triangle inequality yields that for any $t>0$, there holds
\begin{equation*}
 t \| A_q \phi'(t) \|_{L^2(\Omega)} \le \| (tA_q\phi(t))' \|_{L^2(\Omega)} + \| A_q\phi(t)  \|_{L^2(\Omega)} \le c_{T^*} t^{1-\alpha}\| \tilde q-q\|_{L^2\II}.
\end{equation*}
Now the desired inequality follows directly. This completes the proof of the lemma.
\end{proof}

Next we give a stability result. It improves a known result
\cite{ZhangZhou:2017,ZhangZhangZhou:2022} by relaxing the regularity assumption.
\begin{theorem}\label{thm:stab-pot-0}
Let $u_0,f\in  L^2\II$, with $u_0, f \ge m$ a.e. in $\Omega$ and $a_0,a_1\ge m$
for some $m>0$. Then for $q, \tilde q\in \A$, and sufficiently large $T$, there exists $c>0$
independent of $q$, $\tilde q$ and $T$ such that
\begin{equation*}
 \| q - \tilde q \|_{L^2\II} \le c   \|  u(T) - \tilde u(T)  \|_{H^2\II}.
\end{equation*}
\end{theorem}
\begin{proof}
It follows from \eqref{eqn:fde-1} that $q$ can be expressed as
\begin{equation}\label{eqn:q-repres}
  q = [{u(T)}]^{-1}(f-\Dal u(T) + \partial_{xx} u(T)).
\end{equation}
Then we split the difference $q - \tilde q$ into
\begin{align*}
 q - \tilde q 
 & = f\frac{\tilde u(T)-u(T)}{u(T)\tilde u(T)}   + \frac{ u(T)\Dal \tilde u(T) -\tilde u(T)\Dal u(T)  }{u(T)\tilde u(T)} + \frac{\tilde u(T)\partial_{xx} u(T)-u(T)\partial_{xx}\tilde u(T)}{u(T)\tilde u(T)} = \sum_{i=1}^3 {\rm I}_i.
\end{align*}
By the maximum principle of time-fractional diffusion \cite{LuchkoYamamoto:2017}, we deduce $u(T) \ge m>0$. This and the standard Sobolev embedding $H^2(\Omega)\hookrightarrow L^\infty(\Omega)$ (for $d=1,2,3)$ imply
\begin{equation*}
  \|{\rm I}_1\|_{L^2\II} \le m^{-2}\| f \|_{L^2(\Omega)}\| \tilde u(T)-u(T)\|_{H^2\II}.
\end{equation*}
By Lemma \ref{lem:Dalu} and Sobolev embedding theorem, we have the \textit{a priori} bound
$\|u(T)\|_{L^\infty\II} +  \|\Dal u (T)\|_{L^2\II} \le c_{T^*}$ for $T\geq T^*$.
This and Lemma \ref{lem:stab-0} lead to
\begin{align*}
\|{\rm I}_2\|_{L^2\II}
&\le c \big(\|  u(T)\|_{L^\infty} \| \Dal (\tilde u(T) -  u(T)) \|_{L^2\II}  + \| \Dal u(T) \|_{L^2\II}\|  u(T) - \tilde u(T)\|_{L^\infty\II}\big)\\
&\le c \big( T^{-\alpha} \|q-\tilde q\|_{L^2\II} + \|  u(T) - \tilde u(T)\|_{H^2\II} \big),\\
\|{\rm I}_3\|_{L^2\II}&\le c \big(\|  u(T)\|_{L^\infty} \|\partial_{xx} (\tilde u(T) -  u(T)) \|_{L^2\II} + \| \partial_{xx} u(T) \|_{L^2\II}\|  u(T) - \tilde u(T)\|_{L^\infty\II}\big)\\
&\le c  \|  u(T) - \tilde u(T)\|_{H^2\II}.
\end{align*}
Then for sufficiently large $T$, we have
\begin{equation*}
  \| q- \tilde q\|_{L^2\II} \le  c(1-cT^{-\alpha})^{-1}\|  u(T) - \tilde u(T)  \|_{H^2\II} .
\end{equation*}
This completes the proof of the theorem.
\end{proof}

The next stability estimate is the main result of this section.
\begin{theorem}\label{thm:stab-pot}
Let $u_0,f\in  L^2\II$, with $u_0, f \ge m$ a.e. in $\Omega$ and $a_0,a_1\ge m$ for some $m>0$.
Then for $q, \tilde q \in \A$, and $T_0 \le T \le\tilde T$ with sufficiently large $T_0$, there exists $c$ independent
of $q$, $\tilde q$, $T$ and $\tilde T$ such that
\begin{equation*}
 \|  q - \tilde q\|_{L^2\II} \le C \big( \| u(T)- \tilde u(\tilde T) \|_{H^2\II} + T^{-\alpha-1} |T -\tilde T| \big).
\end{equation*}
\end{theorem}
\begin{proof}
In view of the identity \eqref{eqn:q-repres}, we have the following splitting
\begin{align*}
 q - \tilde q & = \Big(\frac{f-\Dal u(T) + \partial_{xx}u(T)}{u(T)} - \frac{f-\Dal \tilde u(T) + \partial_{xx}\tilde u(T)}{\tilde u(T)} \Big)\\
 &\quad\Big( \frac{f-\Dal \tilde u(T) +\partial_{xx} \tilde u(T)}{\tilde u(T)}  - \frac{f-\Dal \tilde u(\tilde T) + \partial_{xx} \tilde u(\tilde T)}{\tilde u(\tilde T)} \Big) =: {\rm I}_1 +{\rm I}_2.
\end{align*}
Theorem \ref{thm:stab-pot-0} implies the following estimate on ${\rm I}_1$:
\begin{align*}
 \| {\rm I}_1  \|_{L^2\II} &\le  c \big( T^{-\alpha} \| q -\tilde q \|_{L^2\II} + \|  u(T) - \tilde u(T)\|_{H^2\II} \big).
\end{align*}
The triangle inequality and Lemma \ref{lem:Dalu} with $\gamma =2$ lead to
\begin{align*}
  \|  u(T) - \tilde u(T)\|_{H^2\II} &\leq \|  u(T) - \tilde u(\tilde T)\|_{H^2\II} +  \|  \tilde u(T) - \tilde u(\tilde T)\|_{H^2\II}\\
   &\leq \|  u(T) - \tilde u(\tilde T)\|_{H^2\II} + c T^{-\alpha-1}|\tilde T - T|,
\end{align*}
and consequently,
\begin{align*}
 \| {\rm I}_1  \|_{L^2\II}
    &\le  c \big( T^{-\alpha} \| q -\tilde q \|_{L^2\II} + \|  u(T) - \tilde u(\tilde T)\|_{H^2\II} + c T^{-\alpha-1}|\tilde T - T|\big).
\end{align*}
Next, to bound the term ${\rm I}_2$, we rewrite
\begin{align*}
{\rm I}_2  & = f\frac{\tilde u(\tilde T)-\tilde u(T)}{\tilde u(T)\tilde u(\tilde T)} + \frac{ \tilde u(T)\Dal \tilde u(\tilde T)- \tilde u(\tilde T)\Dal \tilde u(T)}{\tilde u(T)\tilde u(\tilde T)} + \frac{\tilde u(\tilde T)\partial_{xx} \tilde u(T)-\tilde u(T)\partial_{xx} \tilde u(\tilde T)}{\tilde u(T)\tilde u(\tilde T)} = \sum_{i=1}^3 {\rm I}_{2,i}.
 \end{align*}
By the maximum principle \cite{LuchkoYamamoto:2017}, we have $u(T),\tilde u(\tilde T) \ge m>0$.
This and Lemma  \ref{lem:Dalu} lead to
\begin{equation*}
  \|{\rm I}_{2,1}\|_{L^2\II} \le m^{-2}\| f \|_{L^2(\Omega)}\| \tilde u(\tilde T)-\tilde u(T)\|_{H^2\II}
   \le c T^{-\alpha-1}|\tilde T - T|.
\end{equation*}
Meanwhile, by Sobolev embedding theorem and Lemma  \ref{lem:Dalu}, we obtain
\begin{align*}
  \|{\rm I}_{2,2}\|_{L^2\II} &\le c \big( \| \partial_t^\alpha \tilde u(\tilde T) \|_{L^2(\Omega)} \|  \tilde u(T)-\tilde u(\tilde T)\|_{H^2\II}
 + \| \Dal \tilde u(\tilde T)-\Dal \tilde u(T)\|_{L^2\II} \|  \tilde u(\tilde T)\|_{H^2\II} \big)\\
  &\le c \big( \| \partial_t^\alpha \tilde u(\tilde T) \|_{L^2(\Omega)} + \| \tilde u(\tilde T)\|_{H^2\II}  \big)  \|  \tilde u(T) -\tilde u(\tilde T)\|_{H^2\II} \le c T^{-\alpha-1}|\tilde T - T|,\\
  \|{\rm I}_{2,3}\|_{L^2\II} &\le c \big( \| \partial_{xx} \tilde u(\tilde T) \|_{L^2(\Omega)} \|  \tilde u(T)-\tilde u(\tilde T)\|_{H^2\II} + \| \partial_{xx} (\tilde u(\tilde T)- \tilde u(T))  \|_{L^2\II} \|  \tilde u(\tilde T)\|_{H^2\II} \big)\\
  &\le c \big( \| \partial_t^\alpha \tilde u(\tilde T) \|_{L^2(\Omega)} + \| \tilde u(\tilde T)\|_{H^2\II}  \big) \|\tilde u(T)-\tilde u(\tilde T)  \|_{H^2\II}  \le c T^{-\alpha-1}|\tilde T - T|.
 \end{align*}
Combining the preceding estimates yields
\begin{align*}
\|q - \tilde q\|_{L^2\II}  \le   c T^{-\alpha} \| q -\tilde q\|_{L^2\II} + c \|  u(T) - \tilde u(\tilde T)\|_{H^2\II} + c T^{-\alpha-1}|\tilde T - T|.
\end{align*}
By choosing $T_0$ large enough such that $c T_0^{-\alpha} \le \frac12$, we deduce that for $\tilde T\ge T \ge T_0$, the desired estimate holds.
\end{proof}

The next corollary of Theorem \ref{thm:stab-pot} bounds the terminal time $T$ for perturbed data.
\begin{corollary}
Suppose that $u_0\in L^2\II \backslash \dH {s_1}$,
$f\in \dH s$ and $q\in \A \cap H^s\II$, with $s \in (0, \frac12)$ and $s_1 \in (0, s)$.
Let  $(T, q), \, (\tilde T, \tilde q) \in \mathbb{R}_+\times \A\cap \dH s$ be the solutions of
\texttt{IPP} with observations $u(T)$ and $\tilde u(\tilde T)$, respectively.
Then the following estimate holds
\begin{equation*}
 |T - \tilde T| \le  \Gamma(1-\alpha)^{-\frac1\alpha}\min(\Lambda, \tilde \Lambda)^{-\frac1\alpha-1}|\Lambda - \tilde\Lambda|,
\end{equation*}
with the scalars $\Lambda$ and $ \tilde\Lambda$ respectively given by
\begin{equation*}
\Lambda =\lim_{n\in \mathbb{K} , \ n\rightarrow \infty}
\frac{- (\partial_{xx}u(T), \sin(n\pi x))  }{ (u_0-\phi_0 ,  \sin(n\pi x))}
 \quad
\mbox{and}\quad  \tilde\Lambda = \lim_{n\in \mathbb{K}, \ n\rightarrow \infty}
 \frac{- (\partial_{xx}\tilde u(\tilde T)  , \sin(n\pi x))}{(u_0 - \phi_0, \sin(n\pi x) )} .
\end{equation*}
In particular, for $\Lambda <  \tilde \Lambda$, there holds
$\|q  - \tilde q\|_{L^2\II} \le  c  \big(\| u(T)- \tilde u(\tilde T) \|_{H^2\II} + |\Lambda - \tilde\Lambda|\big)$.
\end{corollary}

\section{Numerical experiments and discussions}\label{sec:numer}

In this section we present numerical results to illustrate simultaneous
recovery of the spatially dependent parameter and terminal time $T$.

\subsection{Numerical algorithm}
First we describe a numerical algorithm for reconstruction. The inverse problems
involve two parameters: unknown time $T$ and space-dependent parameter $v$
($u_0$, $\psi$ or $q$). Since these two parameters
have different influence on the measured data $u(T)$, standard iterative regularization methods, e.g.,
Landweber method and conjugate gradient method, do not work very well. We employ the
Levenberg-Marquardt method \cite{Levenberg:1944,Marquardt:1963}, which has been shown to
be effective for solving related inverse problems \cite{LiaoWei:2019}. Due to the
ill-posedness of the inverse problems, early stopping is required in order to obtain good reconstructions.

Specifically, we define a nonlinear operator
$F:(v,T)\in L^2(\Omega)\times\mathbb{R}_+ \to u(v)(x,T)\in L^2(\Omega)$, where $u(v)$
solves problem \eqref{eqn:fde}, with the parameter $v$. Let $(v^0,T^0)$
be the initial guess of the unknowns $(v,T)$. Now given
the approximation $(v^k,T^k)$, we find the next approximation $(v^{k+1},T^{k+1})$ by
\begin{equation*}
  (v^{k+1},T^{k+1})=\arg\min J_k(v,T),
\end{equation*}
with the functional $J_k(v,T)$ at the $k$th iteration (based at $(v^k,T^k)$) given by
\begin{align*}
  J_k(v,T) =& \frac{1}{2}\|F(v^k,T^k)-g^\delta+\partial_v F(v^k,T^k)(v-v^k)
  +\partial_TF(v^k,T^k)(T-T^k)\|^2_{L^2(\Omega)}\\
   & + \frac{\gamma^{k}}{2}\|v-v^k\|_{L^2(\Omega)}^2
  + \frac{\mu^{k}}{2}|T-T^k|^2,
\end{align*}
where $\gamma^{k}>0$ and $\mu^{k}>0$ are regularization parameters, and $\partial_v F(v^k,T^k)$
and $\partial_TF(v^k,T^k)$ are the derivatives of the forward map $F$ in  $v$ and $T$, respectively.
We employ two parameters since $v$ and $T$ influence the data $u(T)$ differently. The
parameters $\gamma$ and $\mu$ are often decreased geometrically with $\rho\in(0,1)$: $\gamma^{k+1}=\rho\gamma^k$
and $\mu^{k+1}= \rho \mu^k$. The derivative $\partial_v F(v,T)$ can be evaluate
explicitly. For example, for \texttt{BP}, the (directional) derivative
$w=\partial_v F(v,T)[h]$ (in the direction $h$) satisfies
\begin{equation*}
\left\{ \begin{aligned}
 \partial_t^\alpha w -\Delta w + qw &= 0, && \quad \mbox{in }\Omega\times(0,T), \\
 w &= 0, &&\quad \mbox{on }\partial\Omega\times(0,T),\\
 w(0) &= h, &&\quad \mbox{in }\Omega.
\end{aligned}\right.
\end{equation*}
To approximate the derivative $\partial_T F(v,T)$, we use the finite difference
$\partial_T F(v,T) \approx (\delta T)^{-1}(F(v,T+\delta T)-F(v,T))$,
where $\delta T$ is a small number, fixed at $\delta T=1\times10^{-3}$
below. Due to the quadratic structure of the functional $J_k(v,T)$,
the increments $\delta v^k:=v^{k+1}-v^k$ and $\delta T^k:=T^{k+1}-T^k$ satisfy
\begin{align*}
   \left[\begin{array}{cc}
     J_v^*J_v +\mu^{k+1}I & J_v^*J_T\\
     J_T^*J_v  & J_T^*J_T +\gamma^{k+1}
   \end{array}\right] \left[\begin{array}{c}\delta v^k\\ \delta T^k\end{array}\right] = \left[\begin{array}{c}
     J_v^*(g^\delta-F(v^k,T^k))\\
     J_T^*(g^\delta-F(v^k,T^k))
   \end{array}\right],
\end{align*}
where $\ast$ denotes the adjoint operator, and
$J_v = \partial_v F(v^k,T^k)$ and $J_T = \partial_T F(v^k,T^k)$.

\subsection{Numerical illustrations}
Now we present numerical results for the three problems, and with the
domain $\Omega=(0,1)$ in 1D and $\Omega=(0,1)^2$ in 2D, and the terminal time $T=0.5$.
We discretize problem \eqref{eqn:fde} using the Galerkin finite element method with continuous
piecewise linear functions in space, and L1 approximation in time \cite{JinLazarovZhou:2013,
JinLazarovZhou:2016}. The accuracy of a reconstruction $\hat v$ relative to the exact one
$v^\dag$ is measured by the $L^2(\Omega)$ error $e(\hat v)=\|\hat v-v^\dag\|_{L^2(\Omega)}$.
The residual $r(\hat v)$ of the recovered tuple $(\hat v,\hat T)$ is computed as $r(\hat v)
=\|F(\hat v,\hat T)-g^\delta\|_{L^2(\Omega)}$. The exact data $g^\dag$ is generated on a fine
space-time mesh. The noisy data $g^\delta$ is generated from $g^\dag$ by
$g^\delta(x)=g^\dag + \epsilon \|g^\dag\|_{L^\infty(\Omega)}\xi(x)$, where $\xi(x)$
follows the standard Gaussian noise, and $\epsilon>0$ indicates the noise level.

The first example is \texttt{BP}, with $q\equiv0$.
\begin{example}\label{exam:bwp}
\begin{itemize}
  \item[{\rm(i)}]The source $f=\min(x,1-x)$, and the unknown initial condition $u_0=\sin(\pi x)$.
  \item[{\rm(ii)}] The diffusion coefficient $a= 1+\sin(\pi x)y(1-y)$, the source $f=\min(x,1-x)e^x\sin(2\pi y)$, and the
unknown initial condition $u_0 =\sin(\pi x)\sin(\pi y)$.
\end{itemize}
\end{example}

\begin{figure}[hbt]
  \centering
  \begin{tabular}{ccc}
    \includegraphics[width=.32\textwidth]{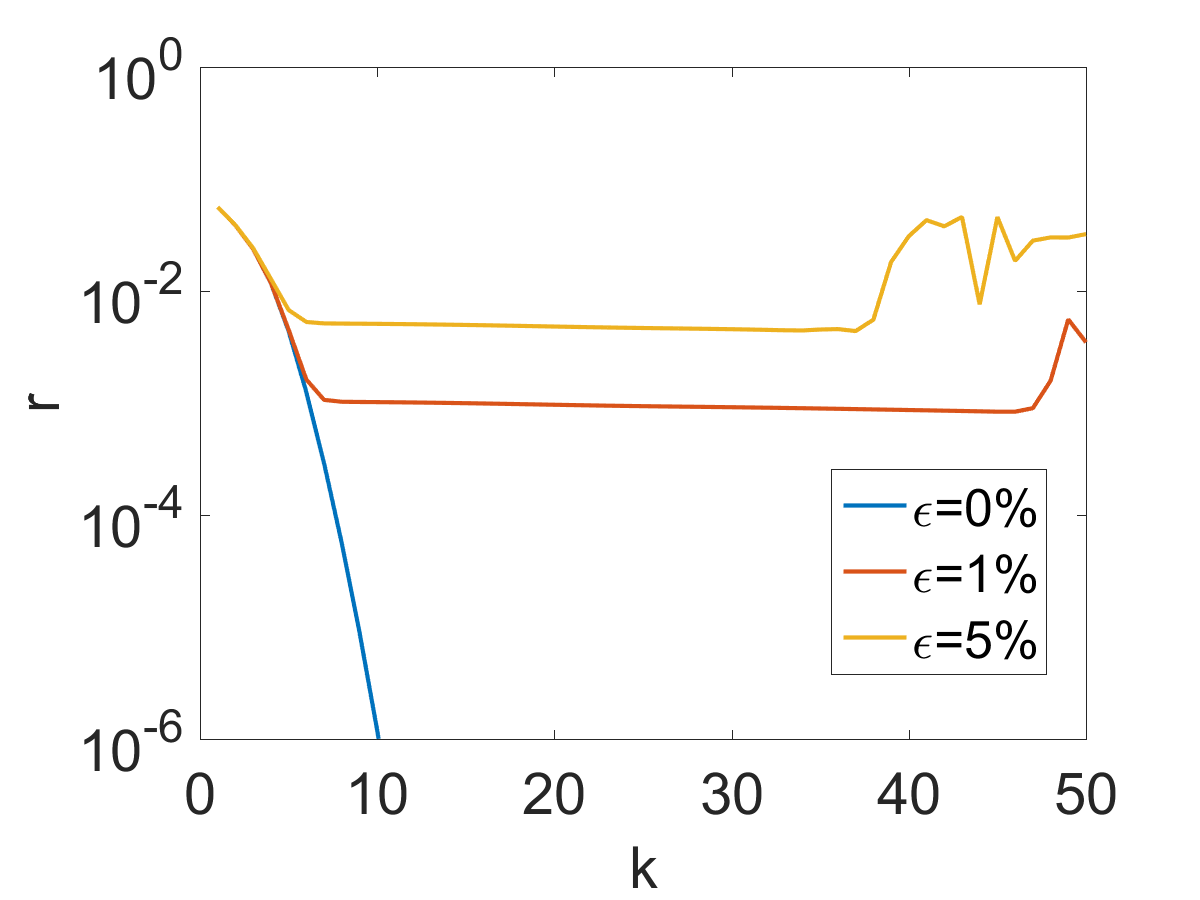}
  & \includegraphics[width=.32\textwidth]{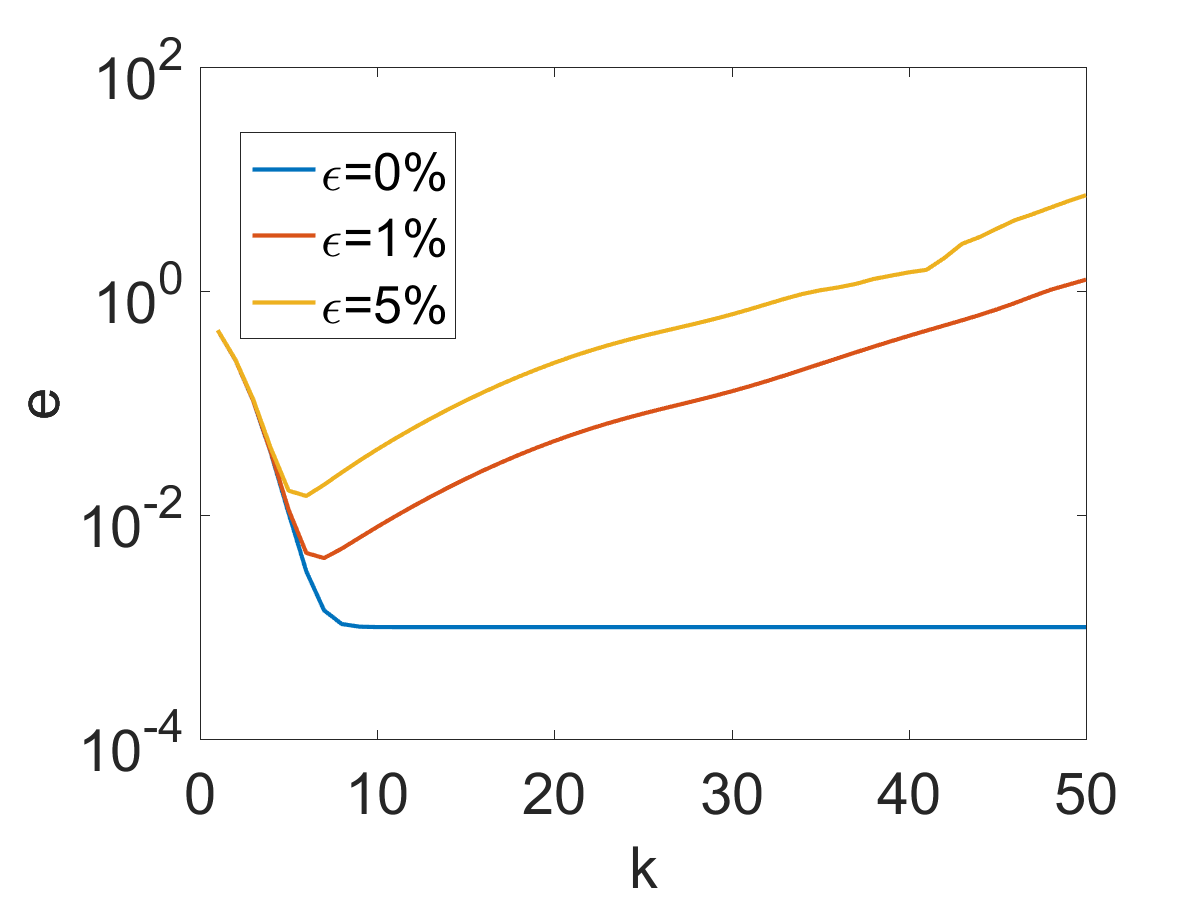}
  & \includegraphics[width=.32\textwidth]{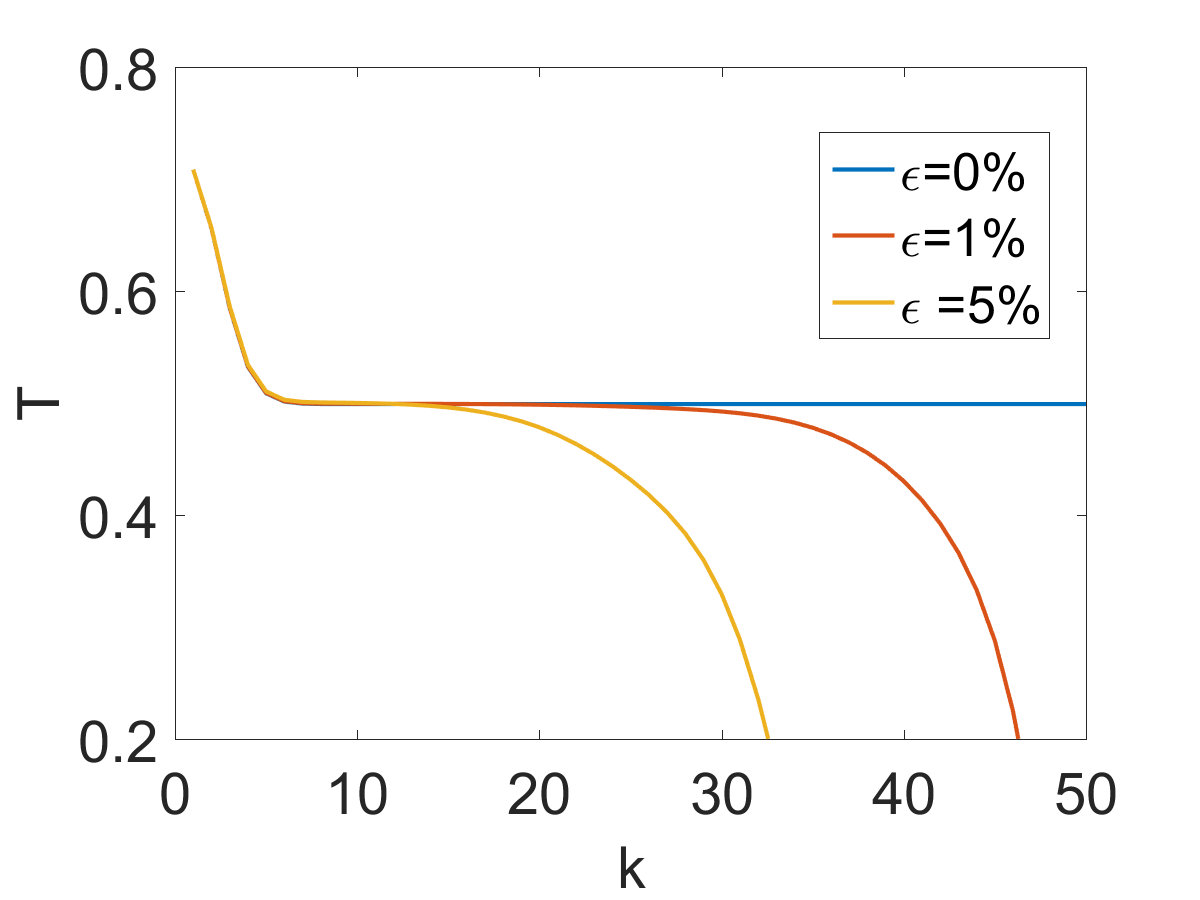}\\
  (a) residual $r$ & (b) error $e$ & (c) terminal time $T$
  \end{tabular}
  \caption{The convergence of Levenberg-Marquadt method for Example \ref{exam:bwp}(i), for $\alpha=0.5$.}\label{fig:bwp-conv}
\end{figure}

\begin{figure}[hbt]
  \centering
  \begin{tabular}{ccc}
    \includegraphics[width=.32\textwidth]{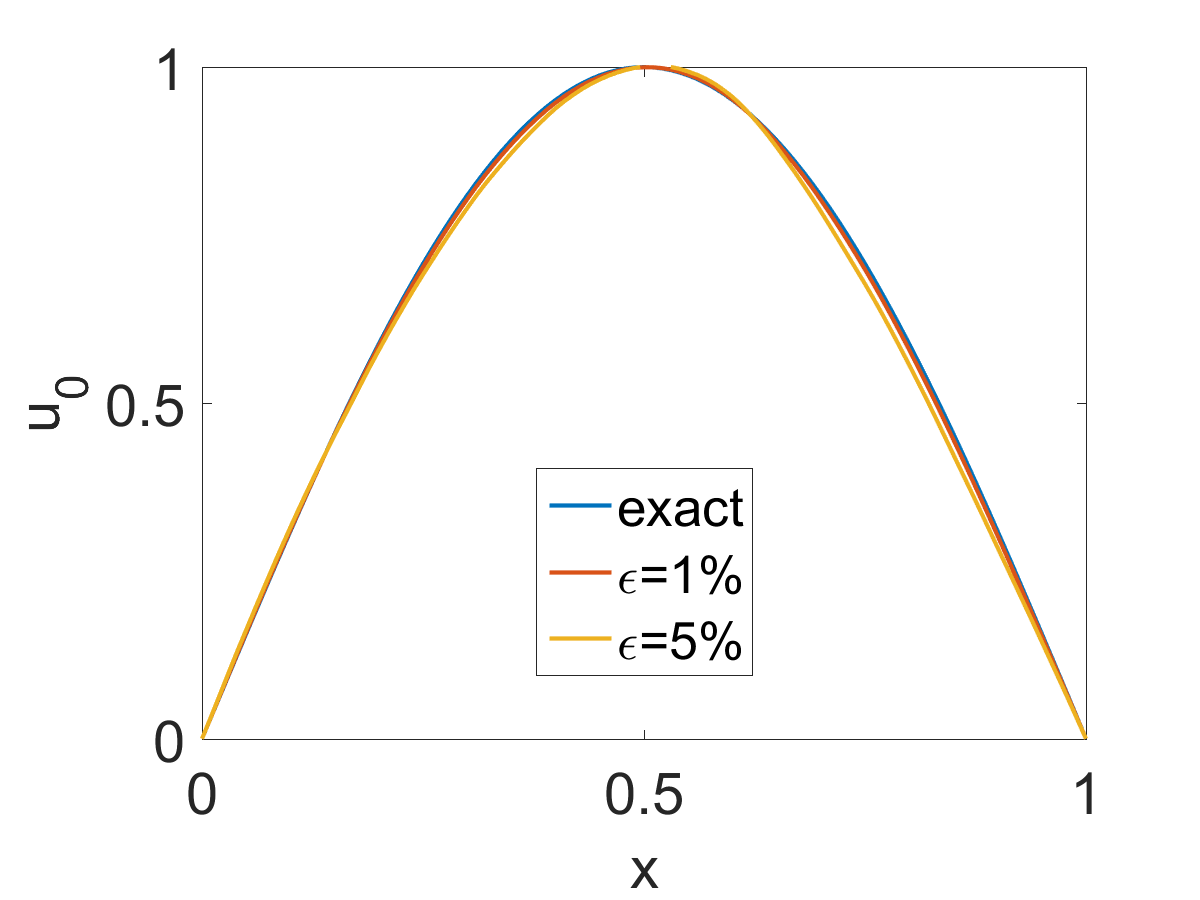}
  & \includegraphics[width=.32\textwidth]{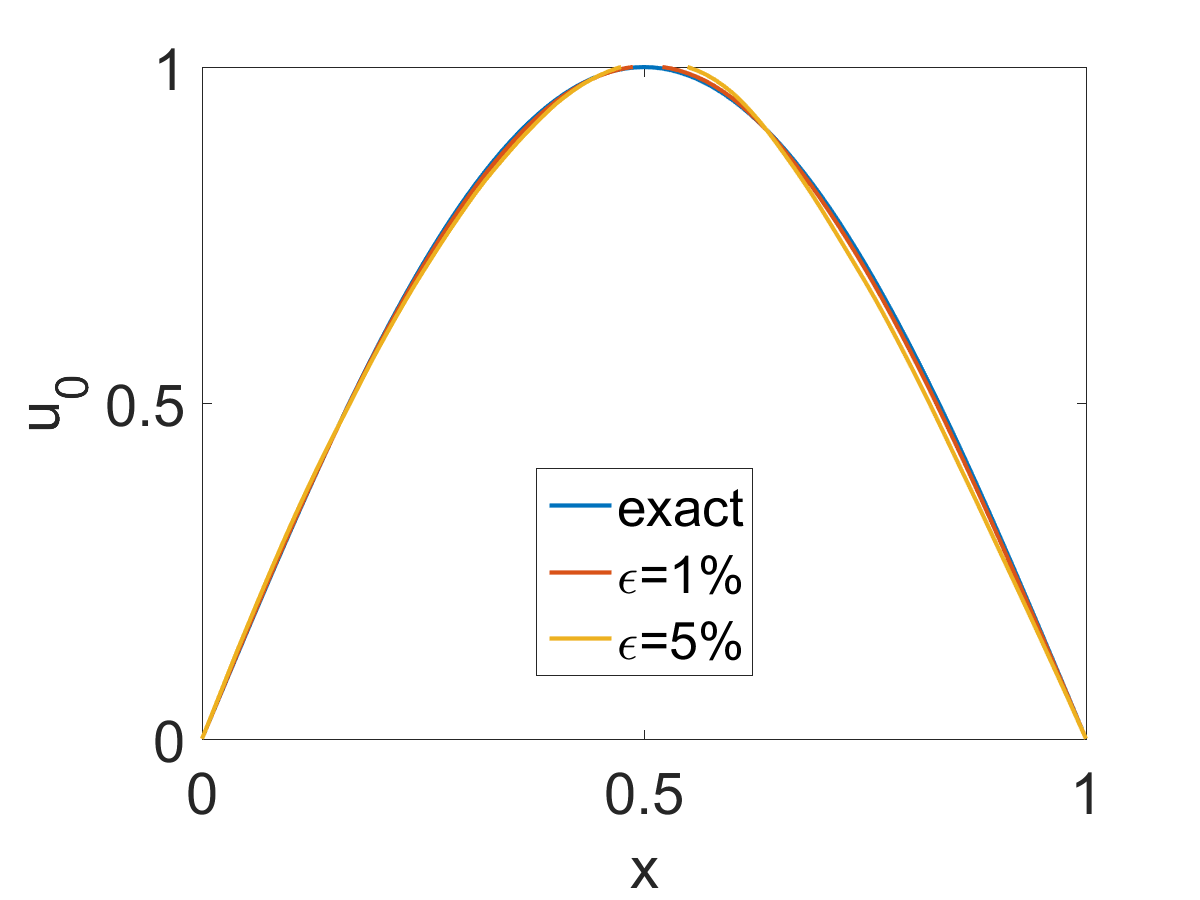}
  & \includegraphics[width=.32\textwidth]{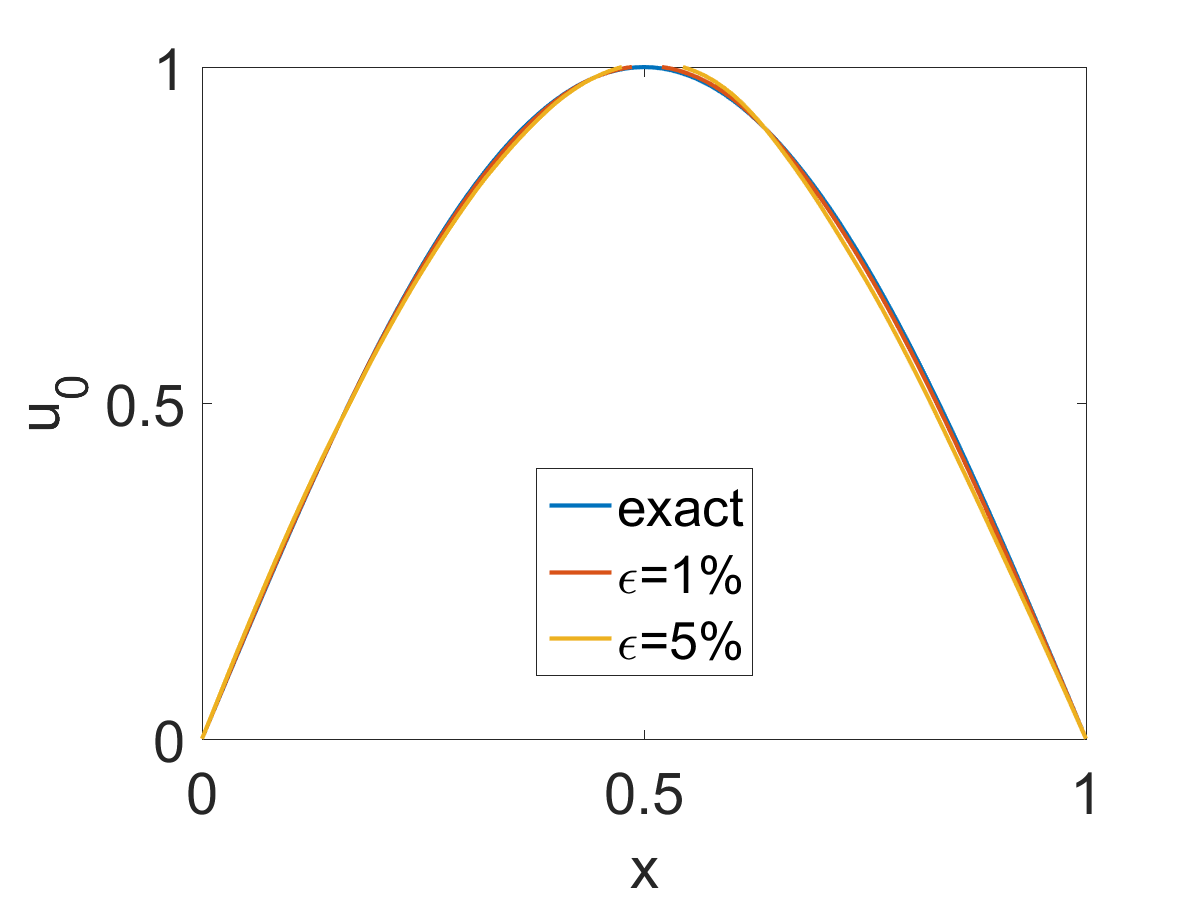}\\
  (a) $\alpha=0.25$ & (b) $\alpha=0.50$ & (c) $\alpha=0.75$
  \end{tabular}
  \caption{The reconstructions of the initial condition $u_0$ for Example \ref{exam:bwp}(i).}\label{fig:bwp-recon}
\end{figure}

\begin{figure}[hbt]
  \centering
  \begin{tabular}{ccc}
    \includegraphics[width=.32\textwidth]{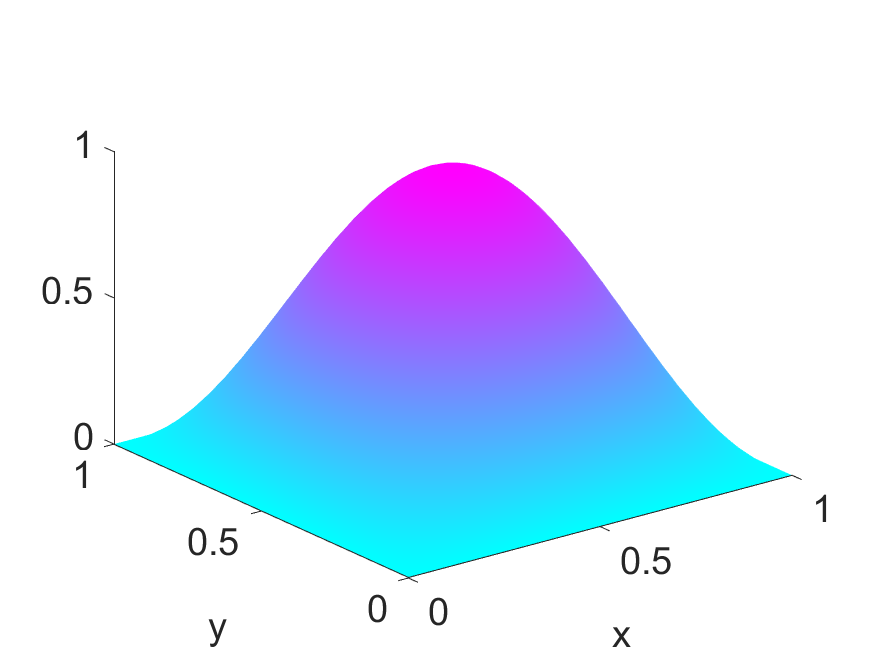}
  & \includegraphics[width=.32\textwidth]{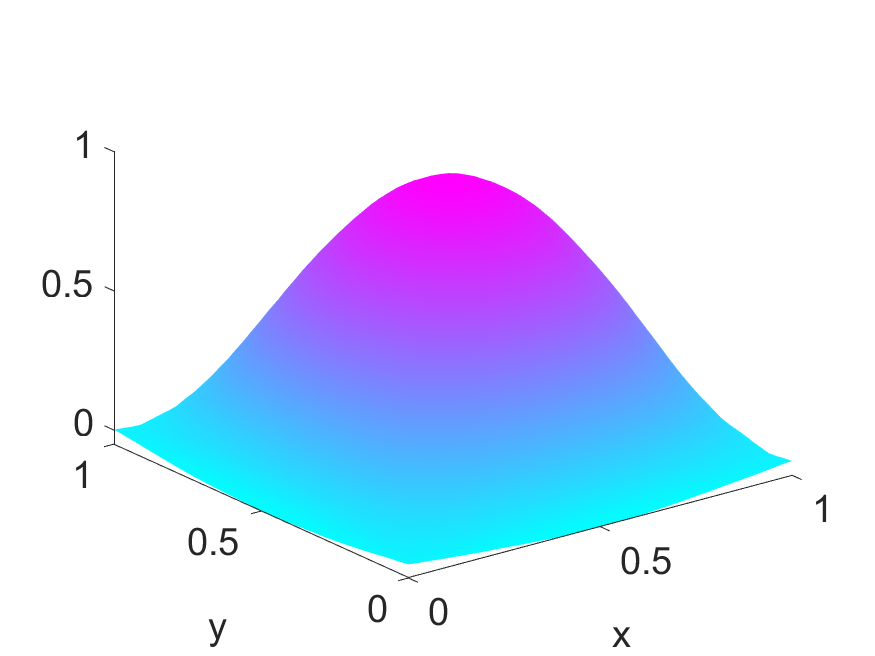}& \includegraphics[width=.32\textwidth]{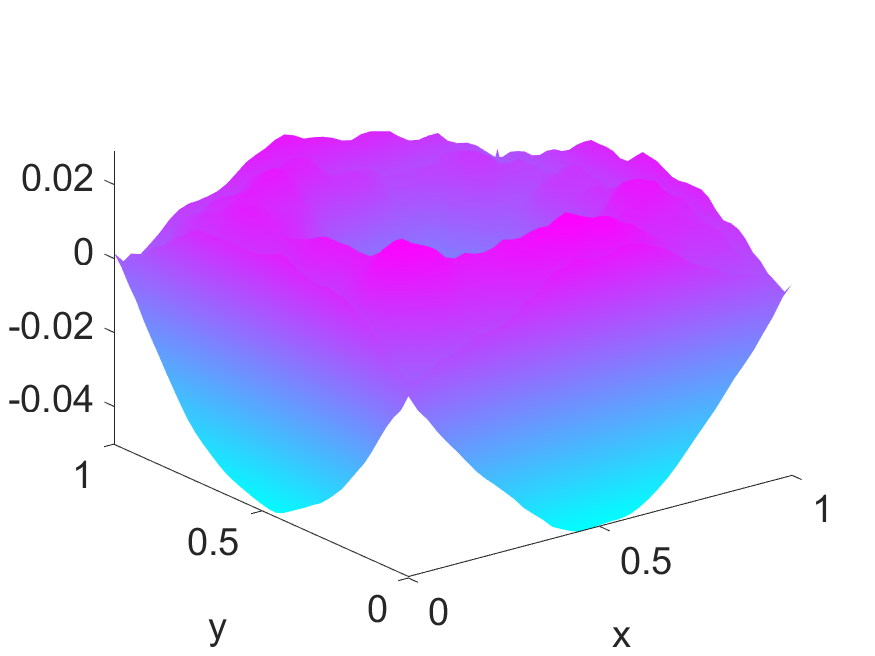}\\
  (a) exact & (b) reconstruction & (c) error
  \end{tabular}
  \caption{The reconstructions of the initial condition $u_0$ for Example \ref{exam:bwp}(ii) with $\alpha=0.5$, $\epsilon=5\%$.}\label{fig:bwp-recon-2d}
\end{figure}

The convergence of the Levernberg-Mardquardt method is shown in Fig. \ref{fig:bwp-conv}.
For the exact data, the residual $r$ decreases rapidly to zero, and the error $e$ also
decreases steadily and eventually levels off at 1e-3 (due to the presence of discretization
error). For noisy data, the method exhibits a typical semi-convergence behavior: the error
$e$ first decreases and then starts to increase rapidly afterwards, necessitating the use of early
stopping. This behavior is also observed for the
estimated terminal time $T$, but it appears to be more resilient to the iteration number $k$ and
it does not change much after a few extra iterations. Nonetheless, the estimate $\hat T$ will
eventually drift away when the method is run for too many iterations. Exemplary reconstructions
of the initial data $u_0$ are shown in Figs. \ref{fig:bwp-recon} and \ref{fig:bwp-recon-2d}, which
has the smallest $L^2(\Omega)$ error along the iteration trajectory; see Table \ref{tab:bwp} for
the stopping index $k^*$. These plots show that the reconstructions are  accurate
for up to $5\%$ noise in the data. See also Table \ref{tab:bwp} for quantitative results. Not that
the accuracy $e$ does not depend very much on the order $\alpha$, and decreases as the noise level
$\epsilon\to0^+$. This agrees with the Lipschitz stability estimate in Theorem \ref{thm:stab-back}.
However, the reconstructions for case (ii) tend to be less accurate than case (i), which is
attributed to discretization errors.

\begin{table}[hbt]
  \centering
  \begin{threeparttable}
  \caption{The numerical results (error $e$, stopping index $k^*$ and recovered $\hat T$) for Example \ref{exam:bwp}.
  For case (i), $\gamma_0=\text{1e-2}$, $\mu_0=\text{2.7e-3}$, $\text{6.3e-3}$ and $\text{1.3e-2}$
for $\alpha=0.25$, $0.50$ and $0.75$, respectively, and $\rho=0.80$; for case (ii), $\gamma_0=\text{1e-3}$,
$\mu_0=\text{1.1e-4}$, $\text{2.7e-4}$, and $\text{6e-4}$ for $\alpha=0.25$, $0.50$ and $0.75$, respectively, and $\rho=0.80$. \label{tab:bwp}}
    \begin{tabular}{ccccccccccc}
    \toprule
    \multicolumn{1}{c}{case} & \multicolumn{1}{c}{$\alpha$}&
    \multicolumn{3}{c}{$0.25$}& \multicolumn{3}{c}{$0.50$}& \multicolumn{3}{c}{$0.75$}\\
    \cmidrule(lr){3-5} \cmidrule(lr){6-8} \cmidrule(lr){9-11}
       &  $\epsilon$   & $e$ &$k^*$ & $\hat T$ &$e$ &$k^*$ & $\hat T$ &$e$ &$k^*$ & $\hat T$ \\
    \midrule
     &   0e-3 & 1.221e-3 &  16 &  0.496 &  9.948e-4 & 20 & 0.498 &  5.023e-4 & 31  &  0.500\\
     &   1e-3 & 1.455e-3 &  8  &  0.496 &  1.210e-3 & 8  & 0.499 &  7.217e-4 & 10  &  0.500\\
(i)  &   5e-3 & 2.961e-3 &  7  &  0.496 &  2.468e-3 & 7  & 0.499 &  1.852e-3 & 9   &  0.500\\
     &   1e-2 & 4.793e-3 &  6  &  0.498 &  4.114e-3 & 7  & 0.499 &  3.275e-3 & 9   &  0.500\\
     &   2e-2 & 8.121e-3 &  6  &  0.498 &  6.872e-3 & 6  & 0.501 &  5.578e-3 & 8   &  0.502\\
     &   5e-2 & 1.675e-2 &  5  &  0.505 &  1.476e-2 & 6  & 0.502 &  1.230e-2 & 8   &  0.503\\
    \midrule\midrule
     &   0e-3 & 7.368e-2  &  32  &  0.496  &  7.414e-2  &  33  &  0.500  &  7.586e-2  &   36  &  0.498\\
     &   1e-3 & 7.574e-2  &  15  &  0.497  &  7.622e-2  &  17  &  0.500  &  7.799e-2  &   20  &  0.498\\
(ii) &   5e-3 & 8.230e-2  &  11  &  0.497  &  8.315e-2  &  13  &  0.500  &  8.648e-2  &   17  &  0.498\\
     &   1e-2 & 9.119e-2  &  10  &  0.497  &  9.247e-2  &  11  &  0.500  &  9.737e-2  &   15  &  0.498\\
     &   3e-2 & 1.169e-1  &   6  &  0.497  &  1.186e-1  &   7  &  0.500  &  1.235e-1  &    7  &  0.499\\
     &   5e-2 & 1.281e-1  &   4  &  0.499  &  1.275e-1  &   4  &  0.503  &  1.300e-1  &    6  &  0.501\\
\bottomrule
\end{tabular}
\end{threeparttable}
\end{table}

The next example is about \texttt{ISP}, with $q\equiv0$.
\begin{example}\label{exam:isp}
\begin{itemize}
  \item[{\rm(i)}] The initial condition $u_0=\sin (2\pi x)$, and the unknown source $\psi(x)=\sin (3\pi x)$.
  \item[{\rm(ii)}] The known diffusion coefficient $a= 1+\sin(\pi x)y(1-y)$, initial condition $u_0 = \sin(\pi x)\sin(\pi y)$, and
the unknown source $\psi= 4x(1-x)e^x\sin(2\pi y)$.
\end{itemize}
\end{example}

\begin{figure}[hbt]
  \centering
  \begin{tabular}{ccc}
    \includegraphics[width=.32\textwidth]{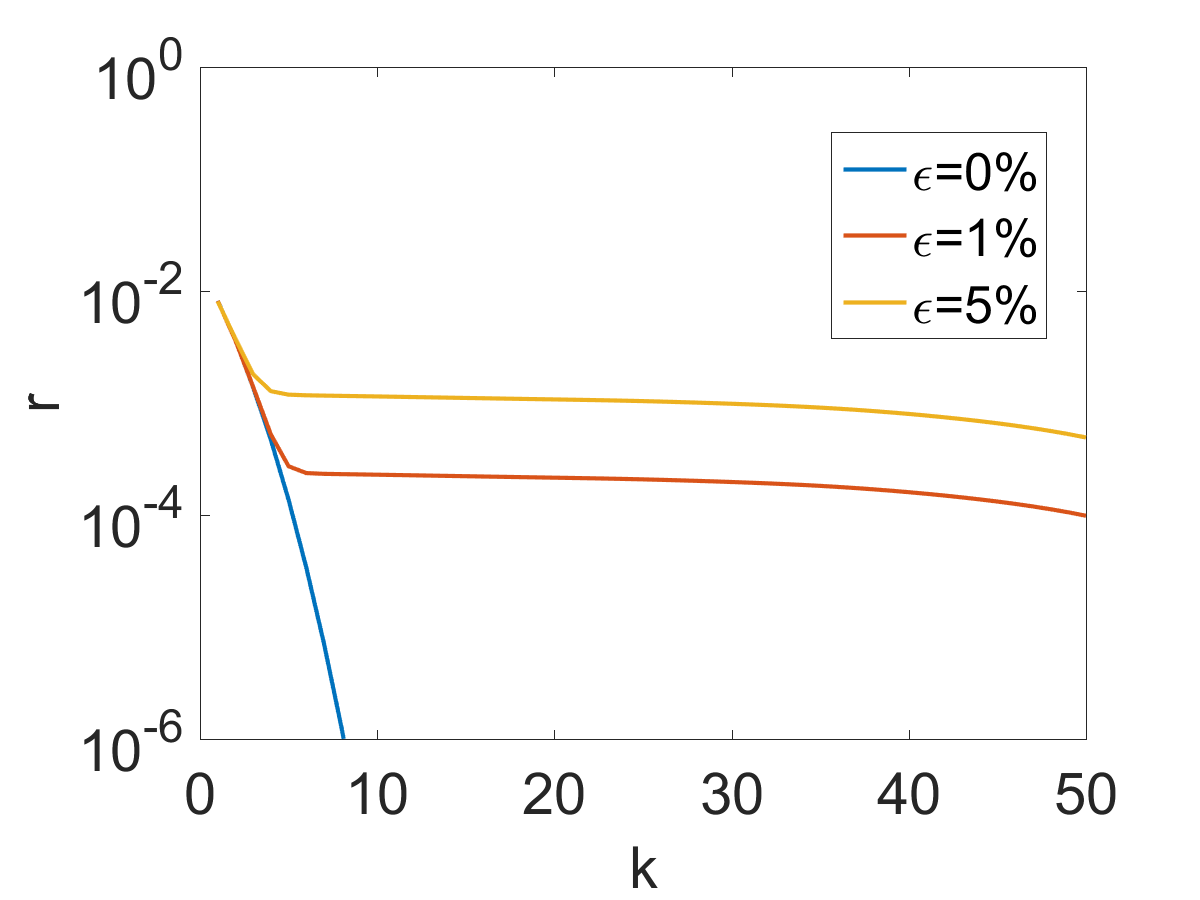}
  & \includegraphics[width=.32\textwidth]{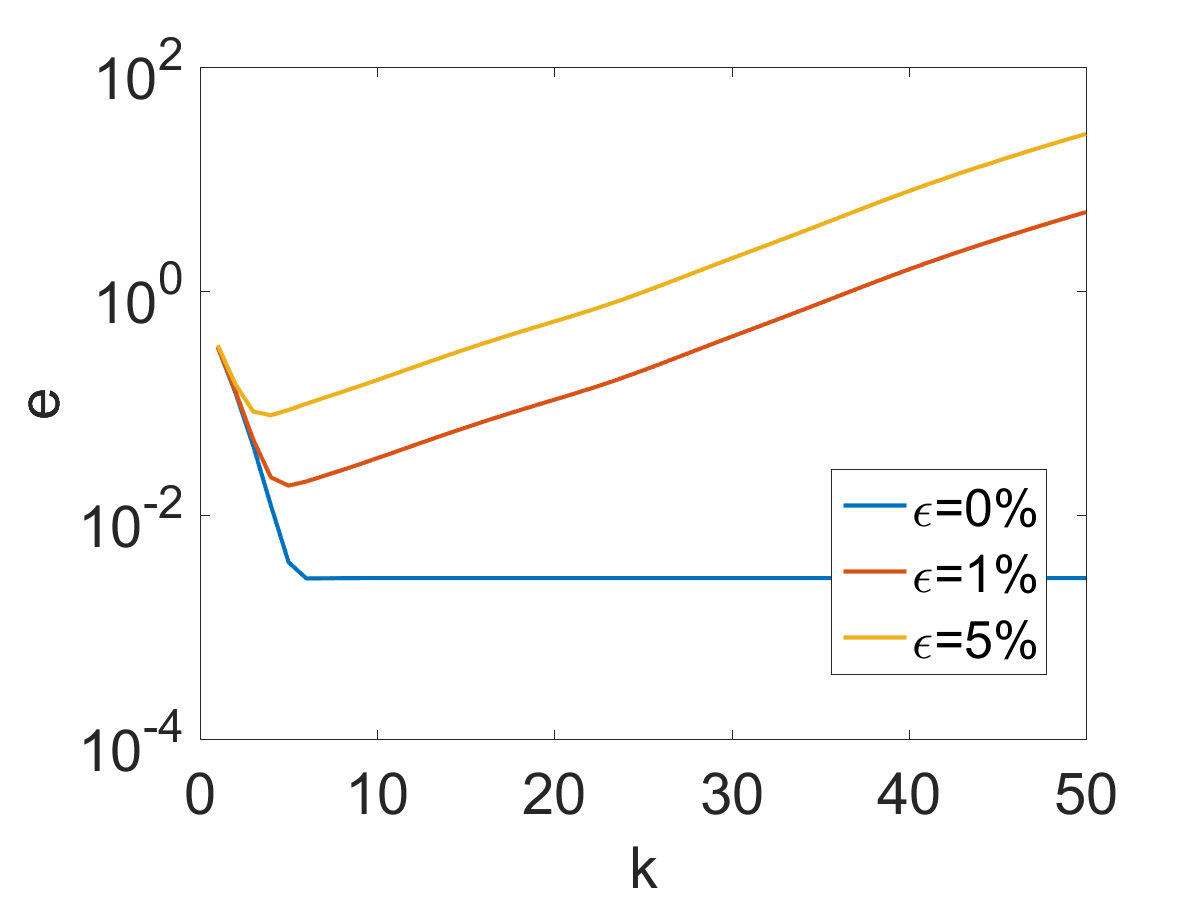}
  & \includegraphics[width=.32\textwidth]{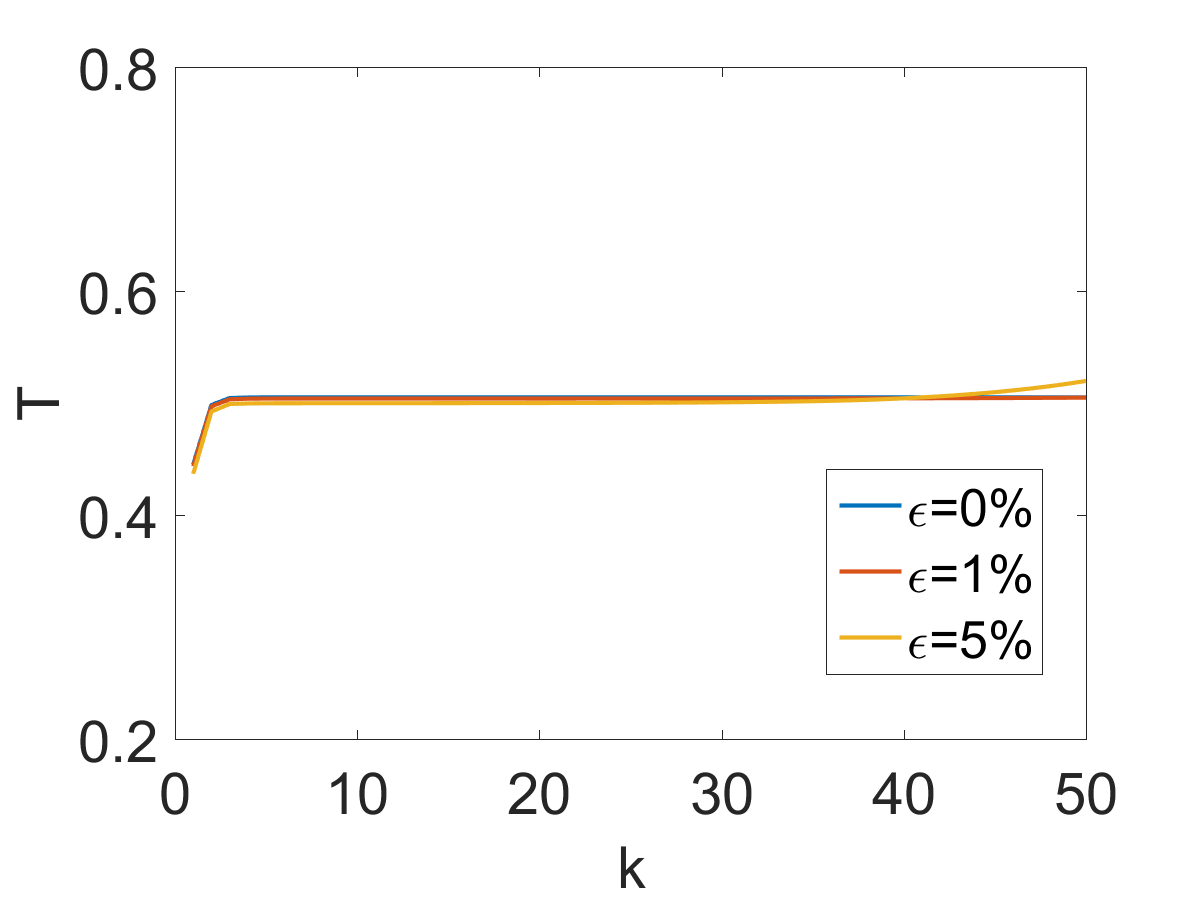}\\
  (a) residual $r$ & (b) error $e$ & (c) terminal time $T$
  \end{tabular}
  \caption{The convergence of the Levenberg-Marquadt method for Example \ref{exam:isp}(i) with $\alpha=0.5$.}\label{fig:isp-conv}
\end{figure}

\begin{figure}[hbt]
  \centering
  \begin{tabular}{ccc}
    \includegraphics[width=.32\textwidth]{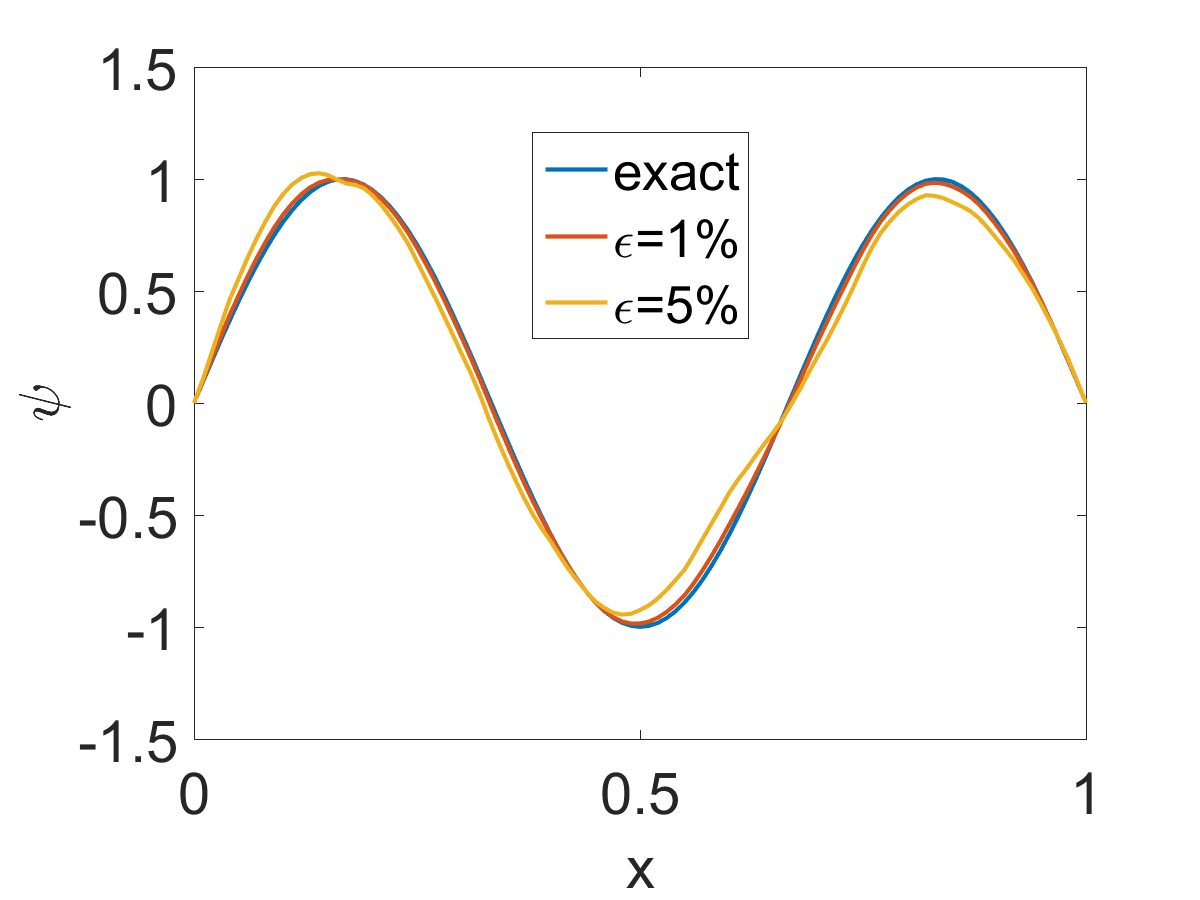}
  & \includegraphics[width=.32\textwidth]{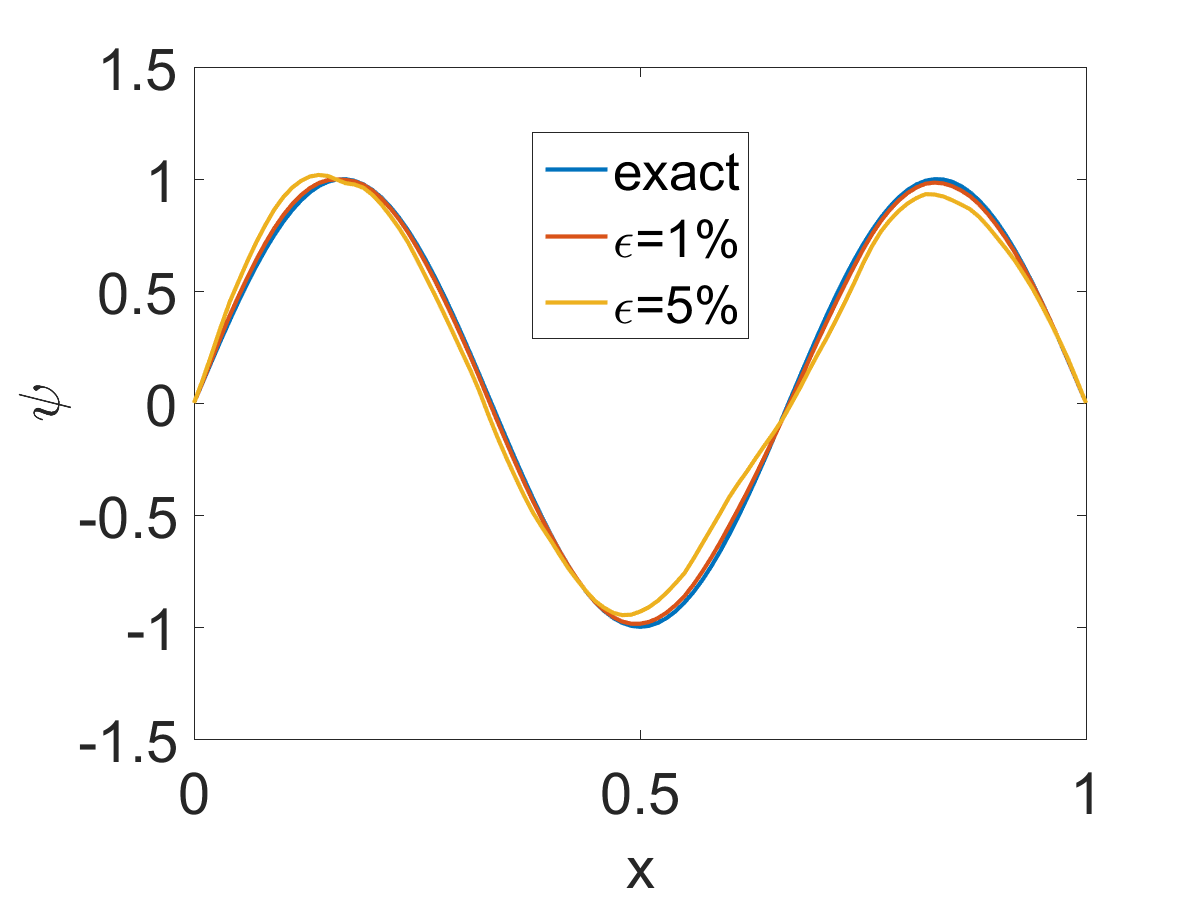}
  & \includegraphics[width=.32\textwidth]{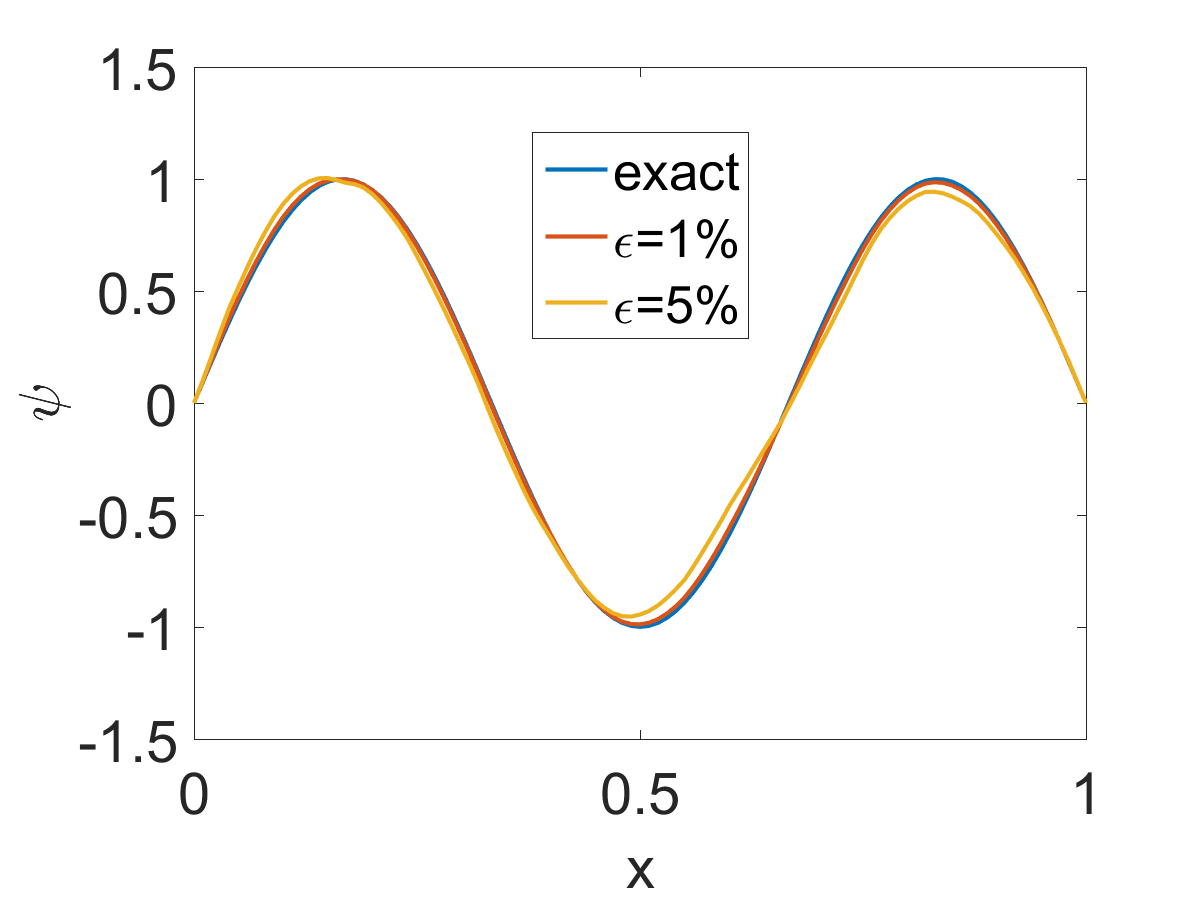}\\
  (a) $\alpha=0.25$ & (b) $\alpha=0.50$ & (c) $\alpha=0.75$
  \end{tabular}
  \caption{The reconstructions of the space-dependent source $\psi$ for Example \ref{exam:isp}(i).}\label{fig:isp-recon}
\end{figure}

\begin{figure}[hbt]
  \centering
  \begin{tabular}{ccc}
    \includegraphics[width=.32\textwidth]{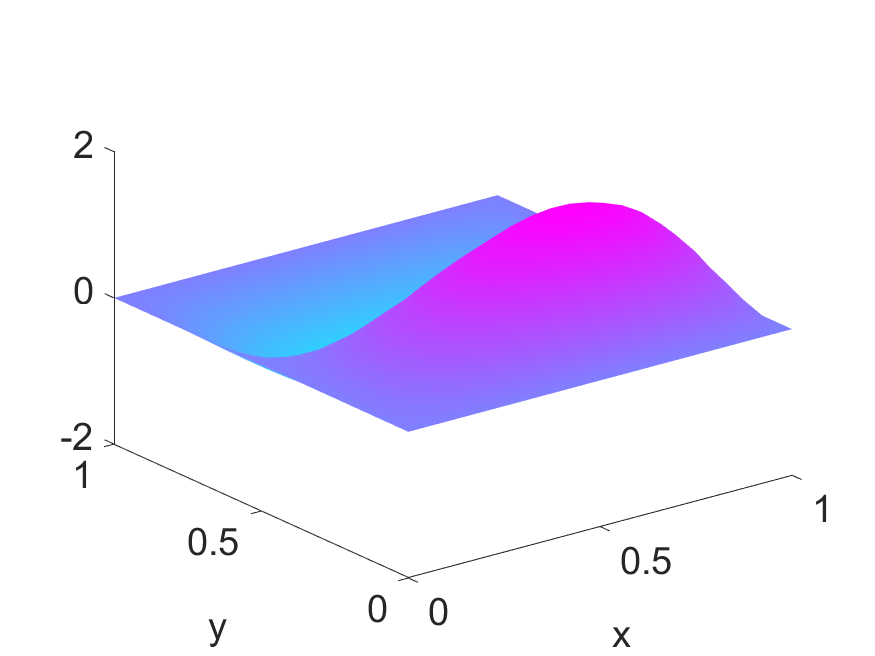}
  & \includegraphics[width=.32\textwidth]{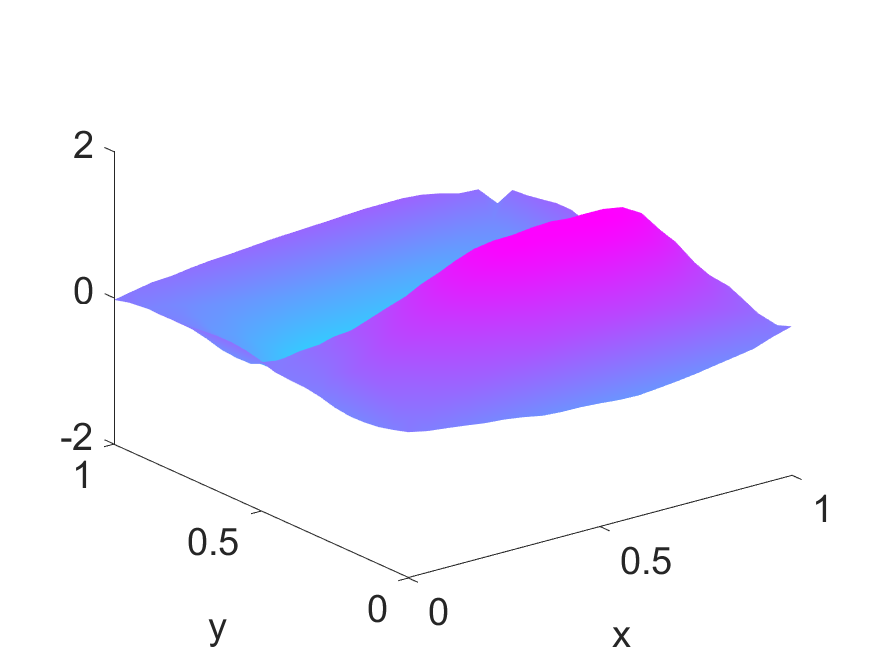} & \includegraphics[width=.32\textwidth]{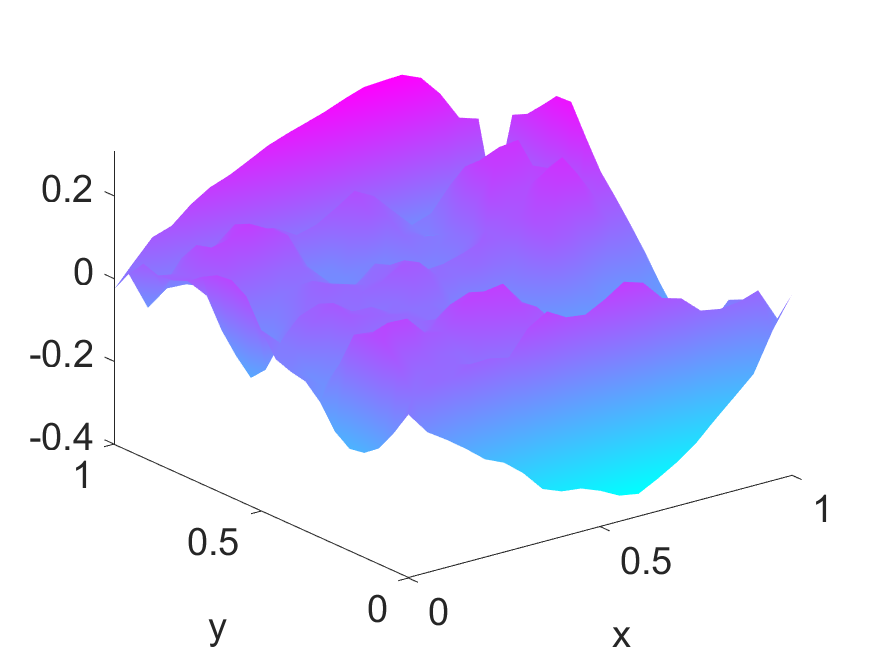}\\
  (a) exact & (b) reconstruction & (c) pointwise error
  \end{tabular}
  \caption{The reconstructions of the space source component $\psi$ for Example \ref{exam:isp}(ii) with $\alpha=0.5$, $\epsilon=5\%$.}\label{fig:isp-recon-2d}
\end{figure}

The numerical results for Example \ref{exam:isp} are shown in Figs. \ref{fig:isp-conv},
\ref{fig:isp-recon} and \ref{fig:isp-recon-2d}, and Table \ref{tab:isp}. The convergence plots
in Fig. \ref{fig:isp-conv} show the semiconvergence phenomenon, and with the chosen parameters,
the method converges rapidly to an acceptable solution, and then the error $e$ starts to increase
shortly afterwards. Nonetheless, the estimate of $T$ converges fairly fast (within 5 iterations),
and it is also quite stable during the iteration. We obtain very accurate reconstructions for
the noise level $\epsilon$ up to $\text{5e-2}$. Like before, the reconstruction quality does not
depend much on the order $\alpha$, cf. Table \ref{tab:isp} and Fig. \ref{fig:isp-recon},
concurring with the observations for \texttt{BP}. The latter also agrees with the fact that
\texttt{ISP} enjoys similar stability as \texttt{BP}, as indicated by Theorems \ref{thm:stab-back}
and \ref{thm:stab-source}.

\begin{table}[hbt]
  \centering
  \begin{threeparttable}
  \caption{The numerical results (error $e$, stopping index $k^*$ and recovered $\hat T$) for Example \ref{exam:isp}.
  For case (i), $\gamma_0=\text{1e-4}$, $\mu_0=\text{1e-8}$, $\text{5e-8}$ and $\text{1e-7}$ for $\alpha=0.25,0.50$
  and $0.75$, respectively, and $\rho=0.8$; and for case (ii), $\gamma_0=\text{1e-4}$, $\mu_0=\text{1e-9}$ and $\rho=0.8$. \label{tab:isp}}
    \begin{tabular}{ccccccccccc}
    \toprule
    \multicolumn{1}{c}{case} & \multicolumn{1}{c}{$\alpha$}&
    \multicolumn{3}{c}{$0.25$}& \multicolumn{3}{c}{$0.50$}& \multicolumn{3}{c}{$0.75$}\\
    \cmidrule(lr){3-5} \cmidrule(lr){6-8} \cmidrule(lr){9-11}
       &  $\epsilon$   & $e$ &$k^*$ & $\hat T$ &$e$ &$k^*$ & $\hat T$ &$e$ &$k^*$ & $\hat T$ \\
    \midrule
     & 0e-3 & 2.534e-3 &  6 &  0.507 &  2.709e-3 & 6 &  0.505 & 3.057e-3 &  6 & 0.506\\
     & 1e-3 & 3.402e-3 &  6 &  0.507 &  3.381e-3 & 6 &  0.504 & 3.417e-3 &  6 & 0.506\\
(i)  & 5e-3 & 1.098e-2 &  5 &  0.506 &  1.000e-2 & 5 &  0.504 & 8.016e-3 &  6 & 0.505\\
     & 1e-2 & 2.035e-2 &  5 &  0.505 &  1.825e-2 & 5 &  0.503 & 1.410e-2 &  5 & 0.505\\
     & 2e-2 & 3.839e-2 &  4 &  0.503 &  3.488e-2 & 4 &  0.502 & 2.669e-2 &  5 & 0.504\\
     & 5e-2 & 8.668e-2 &  4 &  0.496 &  7.761e-2 & 4 &  0.499 & 5.944e-2 &  4 & 0.501\\
\midrule    \midrule
     &   0e-3 & 2.122e-1 & 35 & 0.514 & 2.121e-1 & 35 &  0.508 &2.119e-1 & 35 &  0.509\\
     &   1e-3 & 2.173e-1 &  9 & 0.514 & 2.171e-1 &  9 &  0.508 &2.166e-1 & 10 &  0.509\\
(ii) &   5e-3 & 2.223e-1 &  4 & 0.515 & 2.218e-1 &  4 &  0.509 &2.206e-1 &  5 &  0.509\\
     &   1e-2 & 2.276e-1 &  3 & 0.516 & 2.267e-1 &  3 &  0.509 &2.250e-1 &  3 &  0.509\\
     &   3e-2 & 2.626e-1 &  3 & 0.522 & 2.576e-1 &  3 &  0.512 &2.477e-1 &  3 &  0.512\\
     &   5e-2 & 2.843e-1 &  2 & 0.519 & 2.797e-1 &  2 &  0.507 &2.707e-1 &  2 &  0.500\\
\bottomrule
\end{tabular}
\end{threeparttable}
\end{table}

The last example is about IPP in 1D.
\begin{example}\label{exam:ipp}
The source $f=|\sin(2\pi x)|$, initial condition $u_0\equiv1$, and a zero Dirichlet boundary condition.
The unknown potential $q=\sin^4(\pi x)$.
\end{example}

\begin{figure}[hbt]
  \centering
  \begin{tabular}{ccc}
    \includegraphics[width=.32\textwidth]{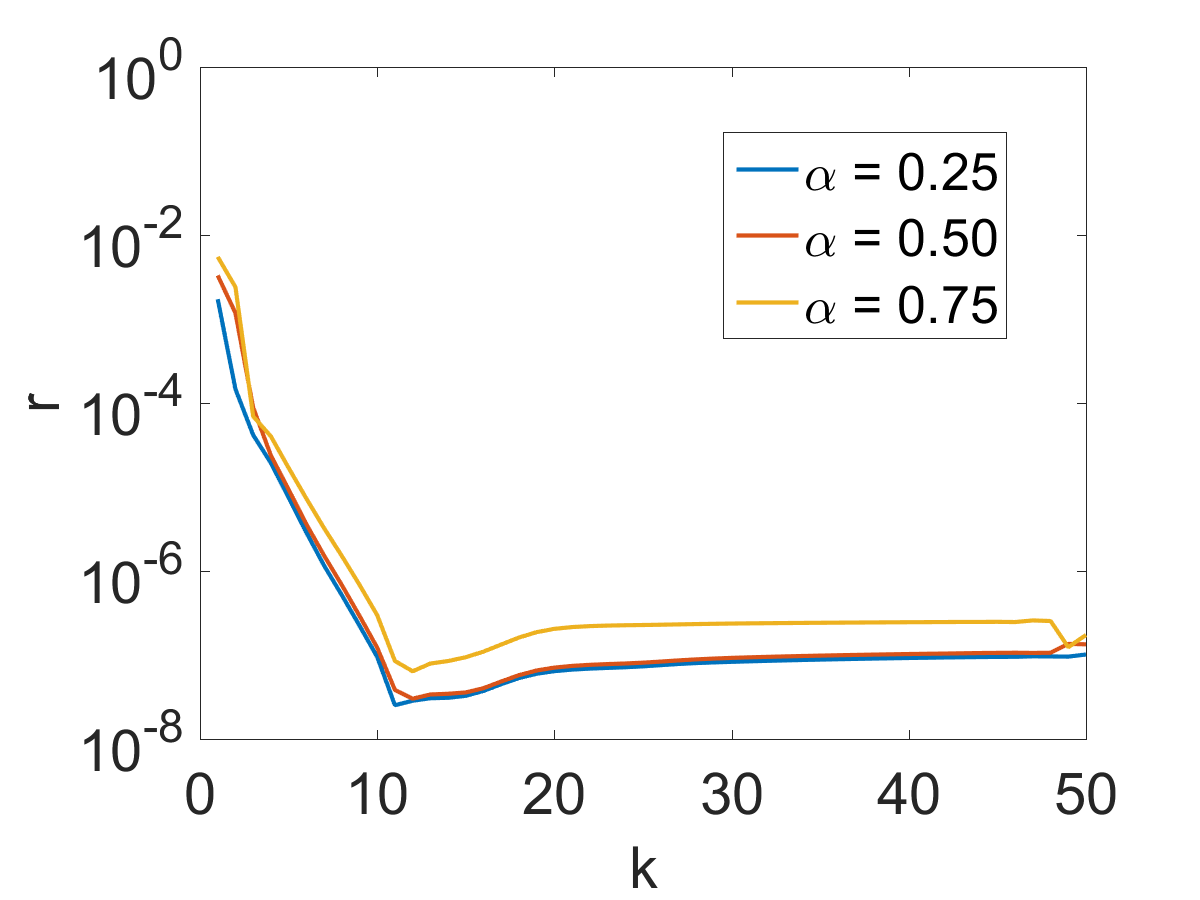}
  & \includegraphics[width=.32\textwidth]{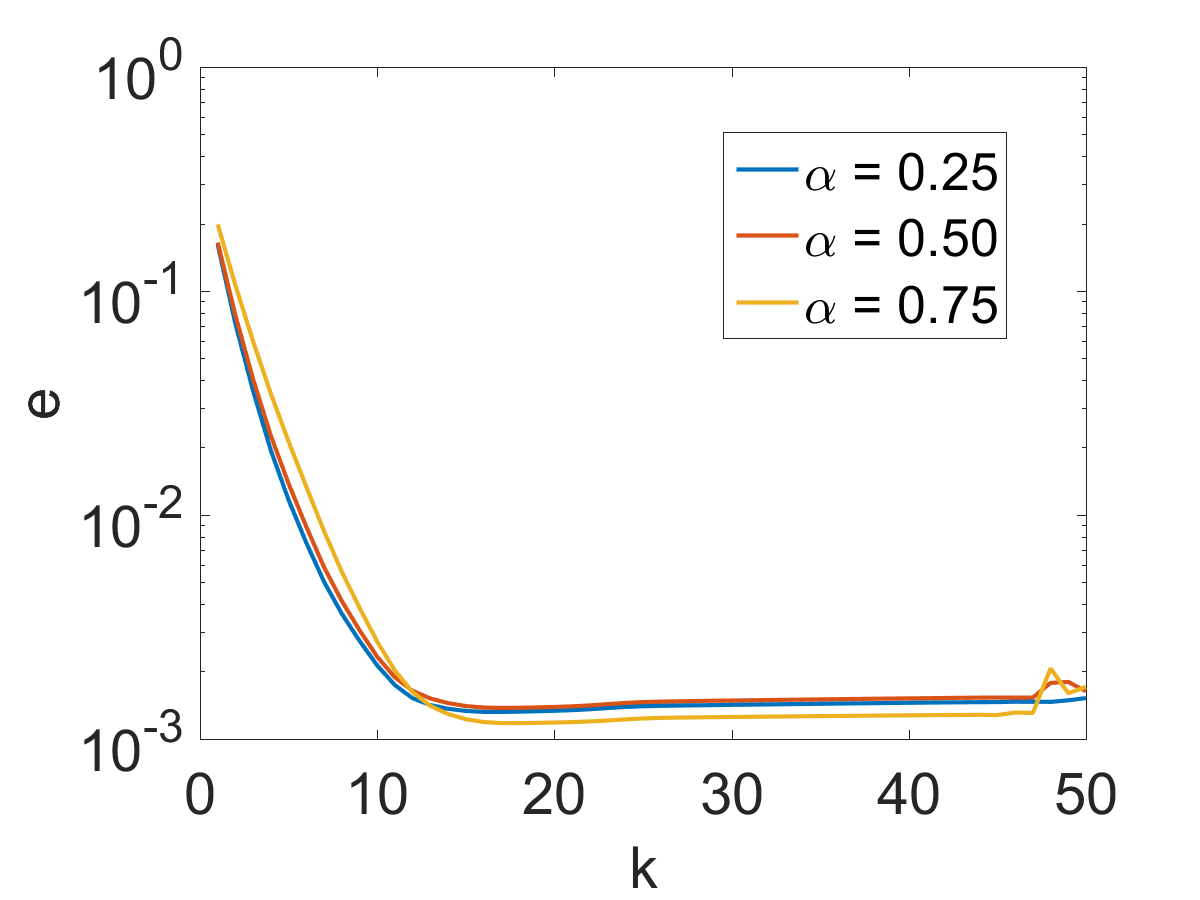}
  & \includegraphics[width=.32\textwidth]{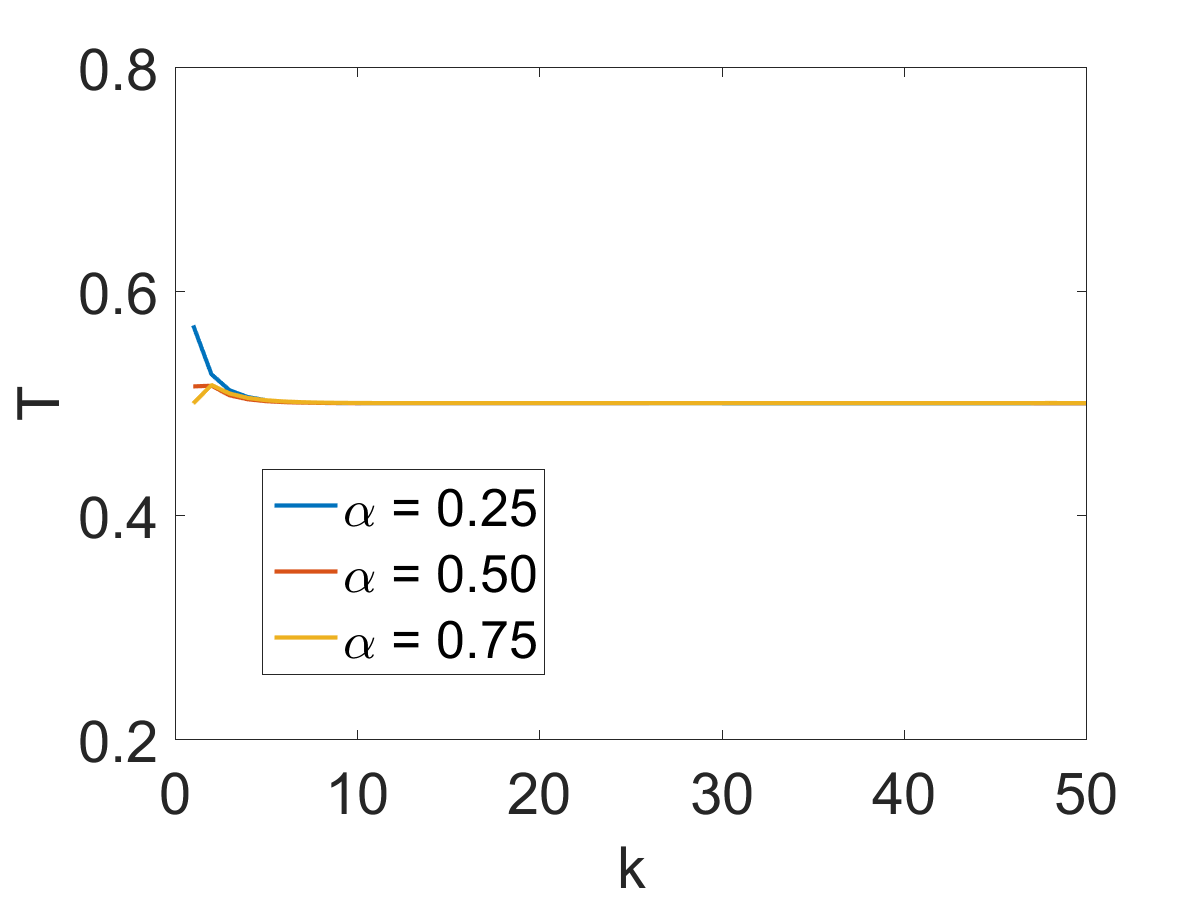}\\
  (a) residual $r$ & (b) error $e$ & (c) terminal time $T$
  \end{tabular}
  \caption{The convergence of the Levenberg-Marquadt method for Example \ref{exam:ipp} with exact data.}\label{fig:ipp-conv}
\end{figure}

\begin{figure}[hbt]
  \centering
  \begin{tabular}{ccc}
    \includegraphics[width=.32\textwidth]{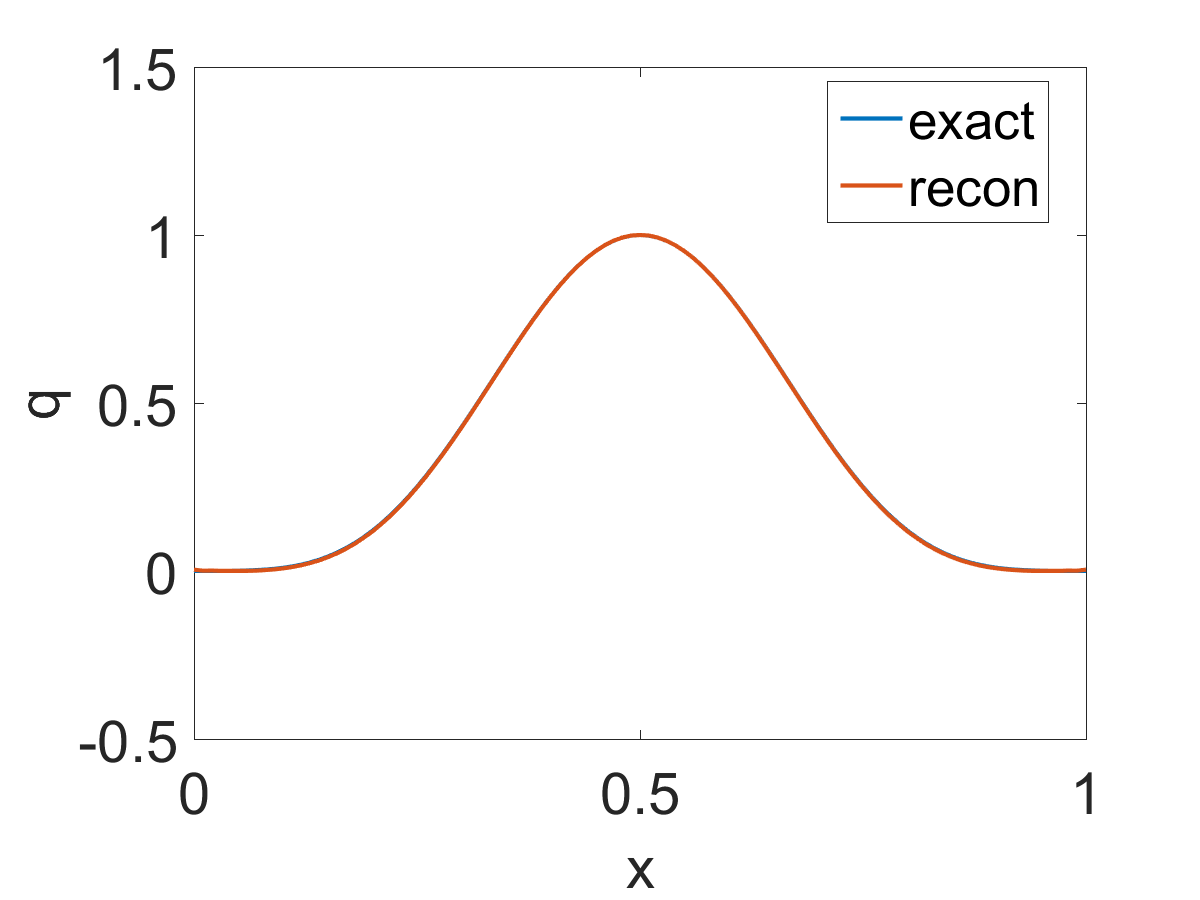}
  & \includegraphics[width=.32\textwidth]{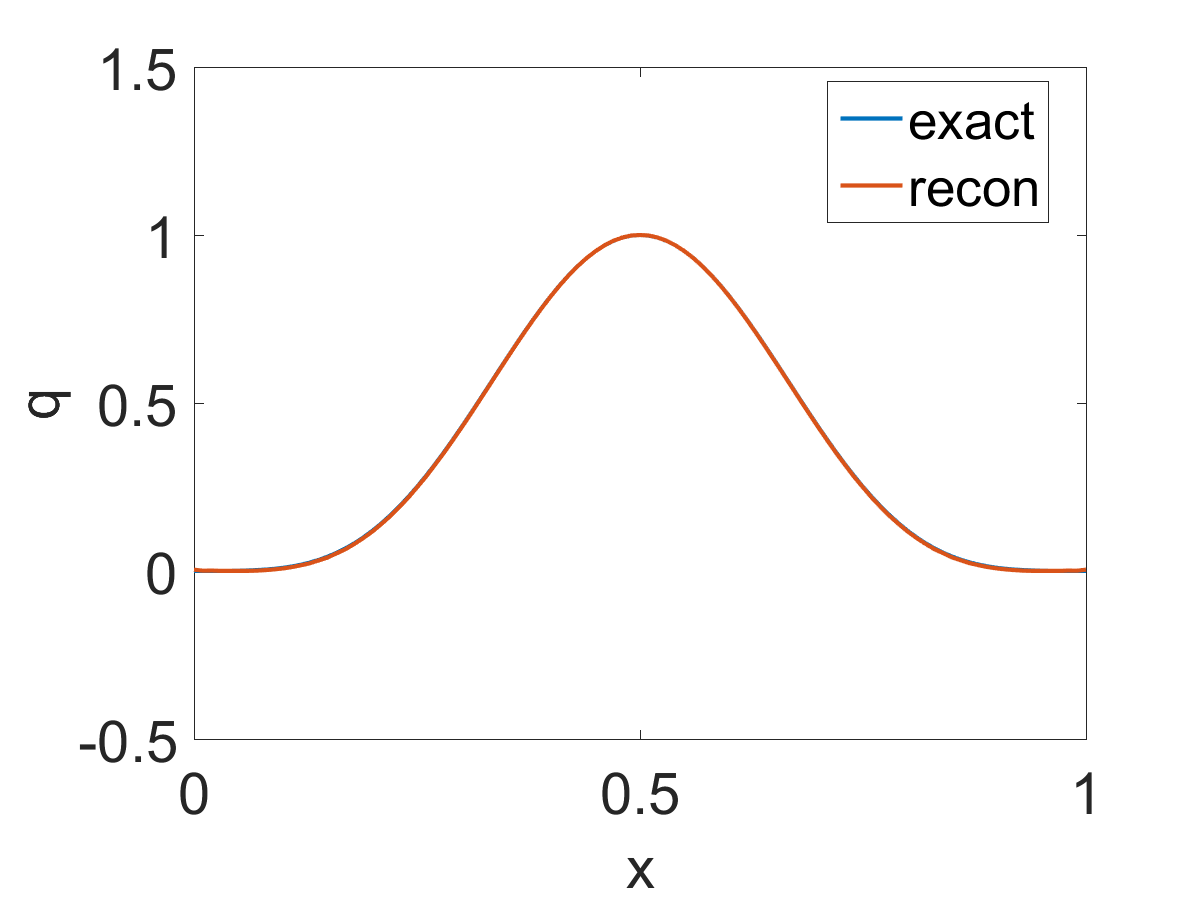}
  & \includegraphics[width=.32\textwidth]{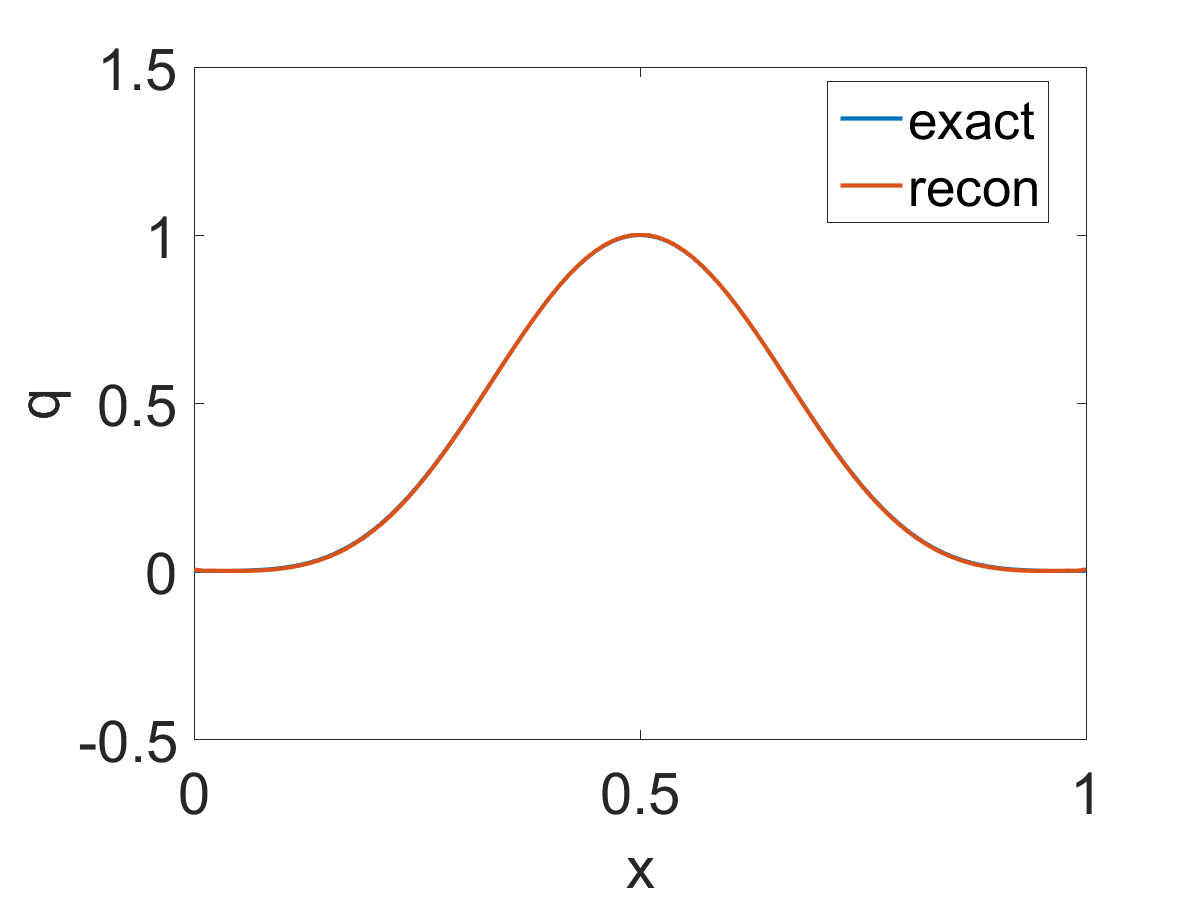}\\
  (a) $\alpha=0.25$ & (b) $\alpha=0.50$ & (c) $\alpha=0.75$
  \end{tabular}
  \caption{The reconstructions of the space-dependent potential $q$ for Example \ref{exam:ipp} with exact data.}\label{fig:ipp-recon}
\end{figure}

The parameter $\gamma_0$ is fixed at $\text{1e-7}$ and $\mu_0$ at $\text{1e-8}$, and the decreasing
factor $\rho$ is set to $0.5$. The numerical results are summarized in Figs. \ref{fig:ipp-conv} and
\ref{fig:ipp-recon}. Note that we can obtain highly accurate reconstructions for exact data, with
the $L^2(\Omega)$ error of the recovered potential $\hat q$ being 1.318e-3, 1.374e-3 and 1.174e-3
for $\alpha=0.25$, 0.50 and 0.75, respectively. The estimated terminal time $T=0.4996$, 0.4997 and
0.4998 for $\alpha=0.25$, 0.50 and 0.75, respectively, are also fairly accurate. This clearly shows
the feasibility of simultaneous recovery. However, for noisy data, the recovery is very challenging.
Numerically we observe that the singular value spectrum of
the linearized forward operator has many tiny values (and hence we have to use a tiny value for the parameter
$\gamma_0$), which precludes applying any realistic amount of noise to the data and renders the recovery from
noisy data highly unstable.

\bibliographystyle{abbrv}
\bibliography{frac}

\begin{thebibliography}{10}

\bibitem{AdamsGelhar:1992}
E.~E. Adams and L.~W. Gelhar.
\newblock Field study of dispersion in a heterogeneous aquifer: 2. spatial
  moments analysis.
\newblock {\em Water Res. Research}, 28(12):3293--3307, 1992.

\bibitem{HatanoHatano:1998}
Y.~Hatano and N.~Hatano.
\newblock Dispersive transport of ions in column experiments: An explanation of
  long-tailed profiles.
\newblock {\em Water Res. Research}, 34(5):1027--1033, 1998.

\bibitem{JannoKinash:2018}
J.~Janno and N.~Kinash.
\newblock Reconstruction of an order of derivative and a source term in a
  fractional diffusion equation from final measurements.
\newblock {\em Inverse Problems}, 34(2):025007, 19, 2018.

\bibitem{Jin:2021book}
B.~Jin.
\newblock {\em {Fractional Differential Equations---An Approach via Fractional
  Dderivatives}}, volume 206 of {\em Applied Mathematical Sciences}.
\newblock Springer, Cham, 2021.

\bibitem{JinLazarovZhou:2013}
B.~Jin, R.~Lazarov, and Z.~Zhou.
\newblock Error estimates for a semidiscrete finite element method for
  fractional order parabolic equations.
\newblock {\em SIAM J. Numer. Anal.}, 51(1):445--466, 2013.

\bibitem{JinLazarovZhou:2016}
B.~Jin, R.~Lazarov, and Z.~Zhou.
\newblock An analysis of the {L}1 scheme for the subdiffusion equation with
  nonsmooth data.
\newblock {\em IMA J. Numer. Anal.}, 36(1):197--221, 2016.

\bibitem{JinRundell:2015}
B.~Jin and W.~Rundell.
\newblock A tutorial on inverse problems for anomalous diffusion processes.
\newblock {\em Inverse Problems}, 31(3):035003, 40, 2015.

\bibitem{JinZhou:2021ip}
B.~Jin and Z.~Zhou.
\newblock An inverse potential problem for subdiffusion: stability and
  reconstruction.
\newblock {\em Inverse Problems}, 37(1):Paper No. 015006, 26, 2021.

\bibitem{KaltenbacherRundell:2019}
B.~Kaltenbacher and W.~Rundell.
\newblock On an inverse potential problem for a fractional reaction-diffusion
  equation.
\newblock {\em Inverse Problems}, 35(6):065004, 31, 2019.

\bibitem{KilbasSrivastavaTrujillo:2006}
A.~A. Kilbas, H.~M. Srivastava, and J.~J. Trujillo.
\newblock {\em Theory and {A}pplications of {F}ractional {D}ifferential
  {E}quations}.
\newblock Elsevier Science B.V., Amsterdam, 2006.

\bibitem{Kou:2008}
S.~C. Kou.
\newblock Stochastic modeling in nanoscale biophysics: subdiffusion within
  proteins.
\newblock {\em Ann. Appl. Stat.}, 2(2):501--535, 2008.

\bibitem{KubicaRyszewskaYamamoto:2020}
A.~Kubica, K.~Ryszewska, and M.~Yamamoto.
\newblock {\em Time-{F}ractional {D}ifferential {E}quations}.
\newblock Springer, Singapore, 2020.
\newblock A theoretical introduction.

\bibitem{Levenberg:1944}
K.~Levenberg.
\newblock A method for the solution of certain non-linear problems in least
  squares.
\newblock {\em Quart. Appl. Math.}, 2:164--168, 1944.

\bibitem{LevitanSargsjan:1975}
B.~M. Levitan and I.~S. Sargsjan.
\newblock {\em Introduction to {S}pectral {T}heory: {S}elfadjoint {O}rdinary
  {D}ifferential {O}perators}.
\newblock American Mathematical Society, Providence, R.I., 1975.

\bibitem{LiYamamoto:2019review}
Z.~Li and M.~Yamamoto.
\newblock Inverse problems of determining coefficients of the fractional
  partial differential equations.
\newblock In {\em Handbook of fractional calculus with applications. {V}ol. 2},
  pages 443--464. De Gruyter, Berlin, 2019.

\bibitem{LiaoWei:2019}
K.~Liao and T.~Wei.
\newblock Identifying a fractional order and a space source term in a
  time-fractional diffusion-wave equation simultaneously.
\newblock {\em Inverse Problems}, 35(11):115002, 23, 2019.

\bibitem{LiuLiYamamoto:2019review-source}
Y.~Liu, Z.~Li, and M.~Yamamoto.
\newblock Inverse problems of determining sources of the fractional partial
  differential equations.
\newblock In {\em Handbook of {F}ractional {C}alculus with {A}pplications.
  {V}ol. 2}, pages 411--429. De Gruyter, Berlin, 2019.

\bibitem{LuchkoYamamoto:2017}
Y.~Luchko and M.~Yamamoto.
\newblock On the maximum principle for a time-fractional diffusion equation.
\newblock {\em Fract. Calc. Appl. Anal.}, 20(5):1131--1145, 2017.

\bibitem{Marquardt:1963}
D.~W. Marquardt.
\newblock An algorithm for least-squares estimation of nonlinear parameters.
\newblock {\em J. Soc. Indust. Appl. Math.}, 11:431--441, 1963.

\bibitem{MillerSamko:1998}
K.~S. Miller and S.~G. Samko.
\newblock A note on the complete monotonicity of the generalized
  {M}ittag-{L}effler function.
\newblock {\em Real Anal. Exchange}, 23(2):753--755, 1997/98.

\bibitem{Nigmatulin:1986}
R.~R. Nigmatulin.
\newblock The realization of the generalized transfer equation in a medium with
  fractal geometry.
\newblock {\em Phys. Stat. Sol. B}, 133:425--430, 1986.

\bibitem{SakamotoYamamoto:2011}
K.~Sakamoto and M.~Yamamoto.
\newblock Initial value/boundary value problems for fractional diffusion-wave
  equations and applications to some inverse problems.
\newblock {\em J. Math. Anal. Appl.}, 382(1):426--447, 2011.

\bibitem{Schneider:1996}
W.~R. Schneider.
\newblock Completely monotone generalized {M}ittag-{L}effler functions.
\newblock {\em Exposition. Math.}, 14(1):3--16, 1996.

\bibitem{Simon:2014}
T.~Simon.
\newblock Comparing {F}r\'{e}chet and positive stable laws.
\newblock {\em Electron. J. Probab.}, 19:no. 16, 25, 2014.

\bibitem{Thomee:2006}
V.~Thom\'{e}e.
\newblock {\em Galerkin {F}inite {E}lement {M}ethods for {P}arabolic
  {P}roblems}.
\newblock Springer-Verlag, Berlin, second edition, 2006.

\bibitem{ZhangZhangZhou:2022}
Z.~Zhang, Z.~Zhang, and Z.~Zhou.
\newblock Identification of potential in diffusion equations from terminal
  observation: analysis and discrete approximation.
\newblock {\em SIAM J. Numer. Anal.}, 60(5):2834--2865, 2022.

\bibitem{ZhangZhou:2017}
Z.~Zhang and Z.~Zhou.
\newblock Recovering the potential term in a fractional diffusion equation.
\newblock {\em IMA J. Appl. Math.}, 82(3):579--600, 2017.

\bibitem{ZhangZhou:2020}
Z.~Zhang and Z.~Zhou.
\newblock Numerical analysis of backward subdiffusion problems.
\newblock {\em Inverse Problems}, 36(10):105006, 27, 2020.

\end{thebibliography}
\end{document}